\documentclass[11pt]{amsart}

\usepackage{amsmath}
\usepackage{amscd}
\usepackage{amssymb}
\usepackage{latexsym}
\usepackage{stmaryrd} 
\usepackage{mathbbol} 

\usepackage[arrow,curve,matrix,tips,2cell]{xy}
  \SelectTips{eu}{10} \UseTips
  \UseAllTwocells
\newtheorem{theorem}{Theorem}[section]
\newtheorem*{theorem*}{Theorem}
\newtheorem{lemma}[theorem]{Lemma}
\newtheorem*{lemma*}{Lemma}
\newtheorem{corollary}[theorem]{Corollary}
\newtheorem*{corollary*}{Corollary}
\newtheorem{proposition}[theorem]{Proposition}

\newtheorem{remark}[theorem]{Remark}
\newtheorem{definition}[theorem]{Definition}

\newcommand{\bgl}{\begin{equation}} 
\newcommand{\egl}{\end{equation}}
\newcommand{\bgloz}{\begin{equation*}} 
\newcommand{\egloz}{\end{equation*}}
\newcommand{\bgln}{\begin{eqnarray}} 
\newcommand{\egln}{\end{eqnarray}}
\newcommand{\bglnoz}{\begin{eqnarray*}} 
\newcommand{\eglnoz}{\end{eqnarray*}}
\newcommand{\btheo}{\begin{theorem}}
\newcommand{\etheo}{\end{theorem}}
\newcommand{\btheooz}{\begin{theorem*}}
\newcommand{\etheooz}{\end{theorem*}}
\newcommand{\blemma}{\begin{lemma}}
\newcommand{\elemma}{\end{lemma}}
\newcommand{\blemmaoz}{\begin{lemma*}}
\newcommand{\elemmaoz}{\end{lemma*}}
\newcommand{\bproof}{\begin{proof}}
\newcommand{\eproof}{\end{proof}}
\newcommand{\bbew}{\begin{beweis}}
\newcommand{\ebew}{\end{beweis}}
\newcommand{\bremark}{\begin{remark}\em}
\newcommand{\eremark}{\end{remark}}
\newcommand{\bdefin}{\begin{definition}}
\newcommand{\edefin}{\end{definition}}
\newcommand{\bprop}{\begin{proposition}}
\newcommand{\eprop}{\end{proposition}}
\newcommand{\bcor}{\begin{corollary}}
\newcommand{\ecor}{\end{corollary}}
\newcommand{\bcoroz}{\begin{corollary*}}
\newcommand{\ecoroz}{\end{corollary*}}
\newcommand{\bfa}{\begin{cases}} 
\newcommand{\efa}{\end{cases}}
\newcommand{\cC}{\mathcal C}

\newcommand{\cF}{\mathcal F}

\newcommand{\cJ}{\mathcal J}

\newcommand{\cL}{\mathcal L}

\newcommand{\cO}{\mathcal O}
\newcommand{\cP}{\mathcal P}

\newcommand{\cR}{\mathcal R}
\newcommand{\cS}{\mathcal S}

\newcommand{\ba}{{\bf a}}
\newcommand{\bb}{{\bf b}}
\newcommand{\bc}{{\bf c}}
\newcommand{\bd}{{\bf d}}
\def\Az{\mathbb{A}}

\def\Nz{\mathbb{N}}

\def\Qz{\mathbb{Q}}

\def\Zz{\mathbb{Z}}

\def\1z{\mathbb{1}}
\newcommand{\fA}{\mathfrak A}

\newcommand{\mfa}{\mathfrak a}
\newcommand{\mfb}{\mathfrak b}
\newcommand{\mfk}{\mathfrak k}
\newcommand{\mfp}{\mathfrak p}
\newcommand{\an}[1]{``#1''} 
\newcommand{\ti}{\tilde}

\newcommand{\lori}{\longrightarrow}

\newcommand{\ma}{\mapsto} 
\newcommand{\mafr}{\mapsfrom} 
\newcommand\onto{\twoheadrightarrow} 
\newcommand\into{\hookrightarrow} 
\newcommand{\Rarr}{\Rightarrow} 
\newcommand{\Larr}{\Leftarrow} 
\newcommand{\LRarr}{\Leftrightarrow} 

\newcommand{\ve}{\varepsilon}

\def\SEMI{\mbox{$\times\kern-2pt\vrule height5pt width.6pt \kern3pt $}}

\newcommand{\halb}{\tfrac{1}{2}}

\newcommand{\img}{{\rm Im\,}}

\newcommand{\Spec}{{\rm Spec\,}} 
\newcommand{\Prim}{{\rm Prim\,}} 
\newcommand{\ind}{{\rm ind\,}}
\newcommand{\res}{{\rm res\,}}
\newcommand{\rk}{{\rm rk\,}}
\newcommand{\ord}{{\rm ord\,}}

\newcommand{\id}{{\rm id}}

\newcommand{\Ad}{{\rm Ad\,}}

\renewcommand{\ker}{{\rm ker}\,}

\newcommand{\reg}{^\times} 
\newcommand{\defeq}{\mathrel{:=}} 

\newcommand{\dop}{\text{: }} 
\newcommand{\falls}{\text{ if }} 
\newcommand{\fa}{\text{ for all }} 
\newcommand{\dotcup}{\ensuremath{\mathaccent\cdot\cup}} 

\newcommand{\Rbar}{\overline{R}} 

\newcommand{\lge}{\left\{} 
\newcommand{\rge}{\right\}} 
\newcommand{\lru}{\left(} 
\newcommand{\rru}{\right)} 
\newcommand{\lsp}{\left\langle} 
\newcommand{\rsp}{\right\rangle} 
\newcommand{\rukl}[1]{\lru #1 \rru} 
\newcommand{\gekl}[1]{\lge #1 \rge} 
\newcommand{\spkl}[1]{\lsp #1 \rsp} 
\newcommand{\menge}[2]{\gekl{ #1 \dop #2 }} 
\newcommand{\ping}{P \subseteq G}

\begin{document}

\title[Semigroup C*-algebras attached to number fields]{On K-theoretic invariants of semigroup C*-algebras attached to number fields}

\author{Xin Li}

\address{Xin Li, Department of Mathematics, Westf{\"a}lische Wilhelms-Universit{\"a}t M{\"u}nster, Einsteinstra{\ss}e 62, 48149 M{\"u}nster, Germany}
\email{xinli.math@uni-muenster.de}


\thanks{\scriptsize{Research supported by the ERC through AdG 267079.}}

\begin{abstract}
We show that semigroup C*-algebras attached to $ax+b$-semigroups over rings of integers determine number fields up to arithmetic equivalence, under the assumption that the number fields have the same number of roots of unity. For finite Galois extensions, this means that the semigroup C*-algebras are isomorphic if and only if the number fields are isomorphic.
\end{abstract}

\maketitle


\setlength{\parindent}{0pt} \setlength{\parskip}{0.5cm}

\section{Introduction}

Recently, there has been some progress in understanding the structure of semigroup C*-algebras (see \cite{Li1}, \cite{Li2}, \cite{Nor}, \cite{C-E-L1}, \cite{C-E-L2}). The construction of these semigroup C*-algebras is easy to describe: Whenever $P$ is a left cancellative semigroup, we can consider the action of $P$ on itself by left multiplication, pass to the induced action of $P$ via isometries on the Hilbert space $\ell^2(P)$ and finally form the C*-algebra of bounded operators on $\ell^2(P)$ generated by all these isometries. This is the reduced semigroup C*-algebra of $P$, denoted by $C^*_r(P)$.

A particular case of this construction, which served as one of the guiding examples for most of the recent work on semigroup C*-algebras, has been the C*-algebra of certain $ax+b$-semigroups: Given a number field $K$ with ring of integers $R$, we take the $ax+b$-semigroup $R \rtimes R\reg$ over $R$ and form the semigroup C*-algebra $C^*_r(R \rtimes R\reg)$. This C*-algebra has been introduced and studied in \cite{C-D-L} (the similar case of the particular semigroup $\Nz \rtimes \Nz\reg$ was treated in \cite{La-Rae} and \cite{La-Nesh}), and it has been further investigated in \cite{Ech-La} and \cite{C-E-L1}. 

Given the natural (and functorial) assignment $K \longmapsto C^*_r(R \rtimes R\reg)$ as described above, the following question immediately comes to mind: What information about the number field $K$ is encoded in the semigroup C*-algebra $C^*_r(R \rtimes R\reg)$? The goal of the present paper is to address this question.

Here are our main results: Let $K$ and $L$ be two number fields, i.e. finite field extensions of $\Qz$. Let $R$ and $S$ be the rings of integers of $K$ and $L$, respectively.
\btheo
\label{main}
Assume that $K$ and $L$ have the same number of roots of unity. If $C^*_r(R \rtimes R\reg) \cong C^*_r(S \rtimes S\reg)$, then $K$ and $L$ are arithmetically equivalent.
\etheo
Arithmetic equivalence for number fields is defined in \cite[\S~1]{Per}.

If we restrict to Galois extensions, we obtain
\btheo
\label{mainGal}
Assume that $K$ and $L$ are finite Galois extensions of $\Qz$ which have the same number of roots of unity. Then we have $C^*_r(R \rtimes R\reg) \cong C^*_r(S \rtimes S\reg)$ if and only if $K \cong L$.
\etheo

The question whether it is possible to read off the number of roots of unity from the semigroup C*-algebra is left open. However, we have the following partial answer:
\btheo
\label{stronger1}
Let $K$ and $L$ be finite Galois extensions of $\Qz$ with rings of integers $R$ and $S$, respectively. Assume that either both $K$ and $L$ have at least one real embedding, or that both $K$ and $L$ are purely imaginary. Then $C^*_r(R \rtimes R\reg) \cong C^*_r(S \rtimes S\reg)$ if and only if $K \cong L$.
\etheo
Alternatively, we could also formulate the following stronger version of Theorem~\ref{mainGal}:
\btheo
\label{stronger2}
Let $K$ and $L$ be finite Galois extensions of $\Qz$ with rings of integers $R$ and $S$, respectively. Assume that either both $K$ and $L$ have only the roots of unity $+1$ and $-1$, or that both $K$ and $L$ have more than two roots of unity. Then $C^*_r(R \rtimes R\reg) \cong C^*_r(S \rtimes S\reg)$ if and only if $K \cong L$.
\etheo

Before we explain the strategy of the proofs, let us first put these results into context. First of all, we should remark that the analogous question for the Bost-Connes system from \cite{Bo-Co}, \cite{Ha-Pa} and \cite{L-L-N} has been studied in \cite{Cor-Mar}, where it is shown that the Bost-Connes system, together with the position of the so-called \an{dagger algebra}, completely determines the number field (up to isomorphism). This is stronger than our result for the semigroup C*-algebra. However, we should point out that the Bost-Connes system is a C*-dynamical system, not only a C*-algebra, and its construction involves not only the number field or its ring of integers as such, but class field theory. 

At the other extreme, it follows from \cite{Cu-Li2}, \cite{Cu-Li3} and \cite{Li-Lu} that if we take the ring C*-algebra of the ring of integers (see \cite{Cu-Li1}, \cite{Cu-Li2}) as such without additional structures, then we obtain very little information about the number field. Since the ring C*-algebra of a ring of integers $R$ is nothing else but the quotient of $C^*_r(R \rtimes R\reg)$ by its maximal primitive ideal, it is certainly not surprising that semigroup C*-algebras contain more information about number fields than ring C*-algebras. The striking observation is that the difference is huge.

Our results could also be interesting from the point of view of classification. First of all, the semigroup C*-algebras we examine here are strongly purely infinite, i.e. $\cO_{\infty}$-absorbing. This follows from \cite[Theorem~8.2.4]{C-E-L1} and \cite[Proposition~2.14]{Pas-Ror}. Therefore, by E. Kirchberg's result in \cite{Kir}, these C*-algebras are classified by KK-theoretic invariants which keep track of the primitive ideal spaces. However, in general it is very complicated to go from these KK-theoretic results to K-theory, i.e, to prove a version of the UCT keeping track of primitive ideal spaces (see \cite{Me-Ne}). From this point of view, it is interesting that for our main results, all we need is ideal related K-theory (with the positions of the units).

For the proof of our main results, the strategy is to extract information about the number field from the K-theory of certain primitive ideals and their quotients of the semigroup C*-algebra. First, it turns out that given a number field $K$ with ring of integers $R$, the minimal non-zero primitive ideals of the semigroup C*-algebra $C^*_r(R \rtimes R\reg)$ correspond one-to-one to the (non-zero) prime ideals of $R$. This is an easy consequence of the results from \cite{Ech-La}. Secondly, we consider the primitive ideal $I_{\gekl{\mfp}}$ corresponding to a prime ideal $\mfp$ of $R$ and form the quotient $C^*_r(R \rtimes R\reg)/I_{\gekl{\mfp}}$. The $K_0$-class of the unit of this quotient (in terms of modules this is the class of the trivial module) turns out to be a torsion element. Building on previous K-theoretic computations in \cite{C-E-L1} and \cite{C-E-L2}, we then show how to read off the prime number $p$ over which $\mfp$ lies (i.e. which satisfies $p \Zz = \mfp \cap \Zz$) from the torsion order of the $K_0$-class of the unit. Once we know this, we can use the result from \cite{dS-Per} that the notions of split and arithmetic equivalence coincide to deduce Theorem~\ref{main} and Theorem~\ref{mainGal}. To deduce the stronger versions Theorem~\ref{stronger1} and Theorem~\ref{stronger2}, we study the maximal primitive ideal of our semigroup C*-algebra $C^*_r(R \rtimes R\reg)$. The corresponding quotient is canonically isomorphic to the ring C*-algebra $\fA[R]$ from \cite{Cu-Li1}. The canonical projection $\pi:C^*_r(R \rtimes R\reg) \onto \fA[R]$ induces a homomorphism in $K_0$, say $\pi_*$. It turns out that for purely imaginary number fields, we can read off the number of roots of unity from the rank of the image of $\pi_*$.

The paper is structured as follows: In a first preliminary section, we gather a few general facts about semigroup C*-algebras and apply these to the particular case of $C^*_r(R \rtimes R\reg)$, $R$ being a ring of integers. In the next section, building on \cite{Ech-La}, we give a very concrete description of the primitive ideals of $C^*_r(R \rtimes R\reg)$. Finally, we turn to K-theoretic invariants attached to the minimal non-zero primitive ideals. In Section~\ref{K-thinv}, we prove Theorem~\ref{main} and Theorem~\ref{mainGal}, and in Section~\ref{number-rou}, we prove Theorem~\ref{stronger1} and Theorem~\ref{stronger2}. The appendix consists of well-known results about restriction homomorphisms in K-theory which are useful for us.

{\bf Acknowledgement:} I am grateful to Gunther Cornelissen for helpful comments about arithmetic equivalence for number fields.

\section{Preliminaries}
\label{pre}

\subsection{Semigroup C*-algebras}

Let $P$ be a left cancellative semigroup with identity element. We first of all recall the construction of the left regular C*-algebra of $P$. Consider the Hilbert space $\ell^2(P)$ with its canonical orthonormal basis $\menge{\varepsilon_x}{x \in P}$, and define for every $p \in P$ an isometry $V_p$ by setting $V_p \varepsilon_x = \varepsilon_{px}$. As in the group case, we simply take the C*-algebra generated by the left regular representation of our semigroup:
\bdefin
$C^*_r(P) \defeq C^* \rukl{\menge{V_p}{p \in P}} \subseteq \cL(\ell^2(P))$.
\edefin

Moreover, it turns out that in the analysis of semigroup C*-algebras as in \cite{Li1}, \cite{Li2}, \cite{C-E-L1} and \cite{C-E-L2}, the following family of right ideals of $P$ plays an important role.
\bdefin
Let $\cJ \defeq \menge{q_1^{-1} p_1 \dotsm q_n^{-1} p_n P}{p_i, q_i \in P} \cup \gekl{\emptyset}$.
\edefin
The right ideals in $\cJ$ are called the constructible right ideals of $P$. A crucial condition on $\cJ$ is given by the independence condition.
\bdefin
\label{ind}
We call $\cJ$ independent if for every $X, X_1, \dotsc, X_n \in \cJ$, the following holds: Whenever $X = \bigcup_{i=1}^n X_i$, then we must have $X = X_i$ for some $1 \leq i \leq n$.
\edefin

Now assume that $P$ is a subsemigroup of a group $G$. In \cite[\S~4]{Li2}, the Toeplitz condition for $\ping$ was introduced. For our purposes, it suffices to note that given a left Ore semigroup $P$, the pair $P \subseteq G = P^{-1} P$ is Toeplitz (see \cite[\S~8.3]{Li2}). Here $G = P^{-1}P$ is the group of left quotients of $P$. In general, given a pair $\ping$ satisfying the Toeplitz condition, we can identify $C^*_r(P)$ with a full corner in a reduced crossed product by $G$ as follows (see \cite[\S~3]{Li2}):

First of all, let $\cJ_{P \subseteq G}$ be the smallest family of subsets of $G$ which contains $\cJ$ and which is closed under left translations by group elements and finite intersections. Next, for every subset $X$ of $G$, let $E_X$ be the orthogonal projection in $\cL(\ell^2(G))$ onto the subspace $\ell^2(X) \subseteq \ell^2(G)$. Alternatively, we can think of $E_X$ as the multiplication operator attached to the characteristic function of $X \subseteq G$. In particular, we can form $E_X$ for all $X \in \cJ_{P \subseteq G}$. We set
\bgloz
  D_{P \subseteq G} \defeq C^*(\menge{E_X}{X \in \cJ_{P \subseteq G}}) \subseteq \ell^{\infty}(G) \subseteq \cL(\ell^2(G)).
\egloz
Here we again think of $\ell^{\infty}(G)$ as acting on $\ell^2(G)$ by multiplication operators.

By construction, the C*-algebra $D_{P \subseteq G}$ is a $G$-invariant sub-C*-algebra of $\ell^{\infty}(G)$, where $G$ acts on $\ell^{\infty}(G)$ by left translations. By \cite[Corollary~3.10]{Li2}, we can identify $C^*_r(P)$ with the full corner $E_P(D_{P \subseteq G} \rtimes_r G)E_P$ of $D_{P \subseteq G} \rtimes_r G$ via
\bgloz
  C^*_r(P) \ni V_p \ma E_P U_p E_P \in E_P(D_{P \subseteq G} \rtimes_r G)E_P.
\egloz
Here, for every $g \in G$, we let $U_g$ be the canonical unitary in the multiplier algebra of $D_{P \subseteq G} \rtimes_r G$ implementing the $G$-action.

Moreover, if $\ping$ is Toeplitz and if $P$ has independent constructible right ideals, \cite[Lemma~4.2]{Li2} tells us that $\cJ_{P \subseteq G}$ is independent (in the sense of \cite[Definition~2.5]{Li2}).

Since $D_{\ping}$ is a commutative C*-algebra, a natural task is to describe its spectrum. This has been done in \cite{Li2}, in the following way: We start with
\bdefin
\label{defSigma}
Let $\Sigma$ be the set of all non-empty $\cJ_{P \subseteq G}$-valued filters, where a $\cJ_{P \subseteq G}$-valued filter is a subset $\cF$ of $\cJ_{P \subseteq G}$ satisfying
\begin{itemize}
\item $X_1 \subseteq X_2 \in \cJ_{\ping}, X_1 \in \cF \Rarr X_2 \in \cF$,
\item $X_1, X_2 \in \cF \Rarr X_1 \cap X_2 \in \cF$,
\item $\emptyset \notin \cF$.
\end{itemize}
\edefin
By \cite[Corollary~2.9]{Li2}, we can identify $\Spec(D_{P \subseteq G})$ and $\Sigma$ via
\bgloz
  \omega: \Spec(D_{P \subseteq G}) \to \Sigma, \chi \ma \menge{X \in \cJ_{P \subseteq G}}{\chi(E_X) = 1}.
\egloz
For this, we used that $\cJ_{P \subseteq G}$ is independent.

Moreover, as explained in \cite{Li2} after Corollary~2.9, the topology of pointwise convergence on $\Spec(D_{P \subseteq G})$ corresponds under $\omega$ to the topology on $\Sigma$ whose basic open sets are given by
\bgloz
  U(X; X_1, \dotsc, X_n) \defeq \menge{\cF \in \Sigma}{X \in \cF, X_i \notin \cF \fa 1 \leq i \leq n}
\egloz
for all $X, X_1, \dotsc, X_n$ in $\cJ_{P \subseteq G}$.

\subsection{The case of $ax+b$-semigroups over rings of integers}

Let $K$ be a number field, i.e. a finite extension of $\Qz$. Let $R$ be the ring of integers in $K$. Moreover, write $K\reg$ for the multiplicative group $K \setminus \gekl{0}$ and $R\reg$ for the multiplicative semigroup $R \setminus \gekl{0}$. We now form the $ax+b$-semigroup $R \rtimes R\reg$. We view it as a subsemigroup of the $ax+b$-group $K \rtimes K\reg$. The main object of interest in this paper is the semigroup C*-algebra $C^*_r(R \rtimes R\reg)$.

We now set out for an explicit description of the dynamical system $G \curvearrowright \Spec(D_{\ping})$ for our pair $P = R \rtimes R\reg \subseteq K \rtimes K\reg = G$. This leads to an explicit identification of $C^*_r(R \rtimes R\reg)$ as a full corner in a crossed product by $K \rtimes K\reg$. Such a description has already been obtained in \cite[\S~5]{C-D-L}. We present a slightly different approach which is more in the spirit of \cite{Li2}. In this way, we set the stage for an explicit description of the primitive ideals of $C^*_r(R \rtimes R\reg)$. This will be the topic of the next section.

We start with the family of constructible right ideals of $R \rtimes R\reg$. It is given by
\bgloz
  \cJ = \menge{(b+\mfa) \times \mfa\reg}{b \in R, (0) \neq \mfa \triangleleft R} \cup \gekl{\emptyset},
\egloz
where $\mfa\reg = \mfa \setminus \gekl{0}$. This is explained in \cite[\S~2.4]{Li1}. Moreover, $\cJ$ is independent by \cite[Lemma~2.30]{Li1}.

Another useful observation is that $R \rtimes R\reg$ is left Ore, and that $K \rtimes K\reg$ is the corresponding group of left quotients. Therefore, our pair $R \rtimes R\reg \subseteq K \rtimes K\reg$ is Toeplitz (see \cite[\S~8.3]{Li2}). The corresponding family $\cJ_{\ping}$ is given by
\bgloz
  \cJ_{P \subseteq G} = \menge{(b+\mfa) \times \mfa\reg}{b \in K, \mfa \text{ a factional ideal of } K} \cup \gekl{\emptyset}.
\egloz

To describe $\Spec(D_{\ping})$ for our pair $R \rtimes R\reg \subseteq K \rtimes K\reg$, consider the finite adele space $\Az_f$ of $K$. The profinite completion $\Rbar$ of $R$ sits as a subring in $\Az_f$. For every $\ba$, $\bb$ in $\Az_f$, we can form the coset $\bb + \ba \Rbar \subseteq \Az_f$.
\bdefin
\label{defcosets}
We set $\cC = \menge{\bb + \ba \Rbar}{\bb, \ba \in \Az_f}$.
\edefin
For $X = (b+\mfa) \times \mfa\reg \in \cJ_{P \subseteq G}$, we let $X(\Rbar)$ be $b + \mfa \cdot \Rbar$, and we set $\emptyset(\Rbar) = \emptyset$. Here we embed $K$ into $\Az_f$ diagonally.

For every prime ideal $\mfp \neq (0)$ of $R$ (in the sequel, by a prime ideal we always mean a non-zero prime ideal of $R$), we can form the discrete valuation ring $R_{\mfp}$ and the corresponding quotient field $K_{\mfp}$ with normalized valuation $v_{\mfp}: K_{\mfp} \to \Zz \cup \gekl{\infty}$. Then
\bgloz
  \Az_f = \menge{(x_{\mfp})_{\mfp} \in \prod_{\mfp} K_{\mfp}}{v_{\mfp}(x_{\mfp}) \geq 0 \text{ for almost all } \mfp}.
\egloz
Every fractional ideal $\mfa$ of $K$ can be factorized in a unique way as $\mfa = \prod_{\mfp} \mfp^{v_{\mfp}(\mfa)}$ with $v_{\mfp}(\mfa) = 0$ for almost all $\mfp$.

\blemma
\label{XR}
For every $X = (b+\mfa) \times \mfa\reg$ in $\cJ_{\ping}$, we have
\bgloz
  X(\Rbar) = b + \mfa \cdot \Rbar = \menge{(x_{\mfp})_{\mfp} \in \Az_f}{v_{\mfp}(-b+x_{\mfp}) \geq v_{\mfp}(\mfa) \fa \mfp}.
\egloz
Moreover, given two sets $X_1$ and $X_2$ in $\cJ_{\ping}$, we have
\bgloz
  X_1(\Rbar) \cap X_2(\Rbar) = (X_1 \cap X_2)(\Rbar).
\egloz
\elemma
\bproof
The first assertion is easy to see. To prove the second claim, let $X_i = (b_i + \mfa_i) \times \mfa_i\reg$ for $i=1,2$. If $X_1 \cap X_2 \neq \emptyset$, then $X_1 \cap X_2 = (b+\mfa) \times \mfa\reg$ with $\mfa = \mfa_1 \cap \mfa_2$ and for some $b \in (b_1 + \mfa_1) \cap (b_2 + \mfa_2)$. Thus
\bglnoz
  && (b_1 + \mfa_1 \cdot \Rbar) \cap (b_2 + \mfa_2 \cdot \Rbar) 
  = (b + \mfa_1 \cdot \Rbar) \cap (b + \mfa_2 \cdot \Rbar) \\
  &=& b + ((\mfa_1 \cdot \Rbar) \cap (\mfa_2 \cdot \Rbar))
  = b + (\mfa_1 \cap \mfa_2) \cdot \Rbar 
  = (X_1 \cap X_2)(\Rbar).
\eglnoz
If $X_1 \cap X_2 = \emptyset$, then $X_1(\Rbar) \cap X_2(\Rbar)$ must be empty as well. If not, then there exists $\bb \in (b_1 + \mfa_1 \cdot \Rbar) \cap (b_2 + \mfa_2 \cdot \Rbar)$. Thus there is $\bb = (b_{\mfp})_{\mfp} \in \Az_f$ with $v_{\mfp}(-b_i + b_{\mfp}) \geq v_{\mfp}(\mfa_i)$ for all prime ideals $\mfp$ and $i=1,2$. But by strong approximation (see \cite[Chapitre~VII, \S~2.4, Proposition~2]{Bour}), there exists $b \in K$ with $v_\mfp(b-b_\mfp) \geq \max(v_\mfp(\mfa_1),v_\mfp(\mfa_2))$ for all prime ideals $\mfp$. It follows that $v_\mfp(-b_i + b) \geq v_\mfp(\mfa_i)$ for all prime ideals $\mfp$ and $i=1,2$. Hence $b$ lies in $(b_1 + \mfa_1) \cap (b_2 + \mfa_2)$.
\eproof

\blemma
We can identify $\Sigma$ and $\cC$ (see Definition~\ref{defSigma} and \ref{defcosets}) via the mutually inverse maps
\bglnoz
  && \Sigma \to \cC, \  \cF \ma \bigcap_{X \in \cF} X(\Rbar) \\
  && \cC \to \Sigma, \  \bb + \ba \Rbar \ma \menge{X \in \cJ_{\ping}}{\bb + \ba \Rbar \subseteq X(\Rbar)}.
\eglnoz
\elemma
\bproof
Let us first prove that these maps are well-defined. For every $\cF \in \Sigma$, the intersection $\bigcap_{X \in \cF} X(\Rbar)$ is non-empty. Namely, if we fix some $Y \in \cF$, then we have
\bgloz
  \bigcap_{X \in \cF} X(\Rbar) = \bigcap_{X \in \cF} X(\Rbar) \cap Y(\Rbar).
\egloz
Now $Y(\Rbar)$ is compact, and $X(\Rbar) \cap Y(\Rbar)$ are closed subsets of $Y(\Rbar)$ satisfying the property that finite intersections are non-empty by the second part of Lemma~\ref{XR} and the definition of a $\cJ_{\ping}$-valued filter. Thus $\bigcap_{X \in \cF} X(\Rbar) \cap Y(\Rbar) \neq \emptyset$. Let us now see why $\bigcap_{X \in \cF} X(\Rbar)$ is of the form $\bb + \ba \Rbar$. As $\bigcap_{X \in \cF} X(\Rbar)$ is non-empty, we can choose $\bb \in \bigcap_{X \in \cF} X(\Rbar)$. For every $X = (b_X + \mfa_X) \times \mfa_X \reg \in \cF$, we have
\bgl
\label{-b+XR}
  - \bb + X(\Rbar) = \mfa_X \cdot \Rbar.
\egl
Define $n_\mfp \defeq \sup_{X \in \cF} v_\mfp(\mfa_X) \in \Zz \cup \gekl{\infty}$ for all every prime ideal $\mfp$. We obviously have $n_\mfp \geq 0$ for almost all $\mfp$. Thus there exists $\ba = (a_\mfp)_\mfp \in \Az_f$ with $v_\mfp(a_\mfp) = n_\mfp$ for all $\mfp$, where we set $v_\mfp(0) = \infty$. This $\ba$ satisfies $\bigcap_{X \in \cF} \mfa_X \cdot \Rbar = \ba \Rbar$. Putting all this together, we obtain
\bgloz
  \bigcap_{X \in \cF} X(\Rbar) = \bb + \bigcap_{X \in \cF} (-\bb + X(\Rbar)) \overset{\eqref{-b+XR}}{=} \bb + \bigcap_{X \in \cF} \mfa_X \cdot \Rbar 
  = \bb + \ba \Rbar.
\egloz
So the first map $\Sigma \to \cC, \  \cF \ma \bigcap_{X \in \cF} X(\Rbar)$ is well-defined. And using the second part of Lemma~\ref{XR}, it is easy to see that $\menge{X \in \cJ_{\ping}}{\bb + \ba \Rbar \subseteq X(\Rbar)}$ is a $\cJ_{\ping}$-valued filter. Thus also the second map $\cC \to \Sigma$, $\bb + \ba \Rbar \ma \menge{X \in \cJ_{\ping}}{\bb + \ba \Rbar \subseteq X(\Rbar)}$ is well-defined. It remains to show that these maps are mutual inverses.

Let us first prove $\cF = \menge{Y \in \cJ_{\ping}}{\bigcap_{X \in \cF} X(\Rbar) \subseteq Y(\Rbar)}$. We obviously have \an{$\subseteq$}. For the reverse inclusion, write every $X \in \cF$ as $X = (b_X + \mfa_X) \times \mfa_X\reg$, and take $Y = (b+\mfa) \times \mfa\reg \in \cJ_{\ping}$ satisfying $\bigcap_{X \in \cF} X(\Rbar) \subseteq Y(\Rbar)$. Then for every prime ideal $\mfp$, we must have $v_\mfp(\mfa) \leq \sup_{X \in \cF} v_\mfp(\mfa_X)$. This implies that for every prime ideal $\mfp$, there exists $X \in \cF$ with $v_\mfp(\mfa) \leq v_\mfp(\mfa_X)$. As $\cF$ is closed under finite intersections, this means that we can find $\ti{X} = (\ti{b} + \ti{\mfa}) \times \ti{\mfa}\reg \in \cF$ with $v_\mfp(\mfa) \leq v_\mfp(\ti{\mfa})$ for all prime ideals $\mfp$, i.e. $\mfa \supseteq \ti{\mfa}$. As $\ti{X}$ and $Y$ satisfy $\bigcap_{X \in \cF} X(\Rbar) \subseteq \ti{X}(\Rbar)$ and $\bigcap_{X \in \cF} X(\Rbar) \subseteq Y(\Rbar)$, the intersection $\ti{X} \cap Y$ cannot be empty by the second part of Lemma~\ref{XR}. As $\ti{\mfa} \subseteq \mfa$, we conclude that $\ti{X}$ must be contained in $Y$, hence $Y \in \cF$. This proves \an{$\supseteq$}.

Finally, let us prove $\bigcap \menge{X \in \cJ_{\ping}}{\bb + \ba \Rbar \subseteq X(\Rbar)} = \bb + \ba \Rbar$. It is clear that \an{$\supseteq$} holds. To prove \an{$\subseteq$}, we use our observation that $\bigcap \menge{X \in \cJ_{\ping}}{\bb + \ba \Rbar \subseteq X(\Rbar)}$ is of the form $\bd + \bc \Rbar$ for some $\bd$, $\bc$ in $\Az_f$. It suffices to prove $\bc \Rbar \subseteq \ba \Rbar$ since we already know $\bd + \bc \Rbar \supseteq \bb + \ba \Rbar$. Let $\cF = \menge{X \subseteq \cJ_{\ping}}{\bb + \ba \Rbar \subseteq X(\Rbar)}$ and write every $X \in \cF$ as $X = (b_X + \mfa_X) \times \mfa_X\reg$. Then we know that $v_\mfp(\bc) = \sup_{X \in \cF} v_\mfp(\mfa_X)$. We also know that $v_\mfp(\mfa_X) \leq v_\mfp(\ba)$ for all prime ideals $\mfp$ as $\bb + \ba \Rbar \subseteq X(\Rbar)$. Hence it follows that if $v_\mfp(\bc) = \infty$, then also $v_\mfp(\ba) = \infty$. If $v_\mfp(\ba)$ is finite, then there exists $X \in \cF$ with $v_\mfp(\mfa_X) = v_\mfp(\ba)$, and we deduce $v_\mfp(\bc) \geq v_\mfp(\ba)$. Thus $\bc \Rbar \subseteq \ba \Rbar$ as desired. This proves \an{$\subseteq$}.
\eproof

Let us denote the map $\Sigma \to \cC$, $\cF \ma \bigcap_{X \in \cF} X(\Rbar)$ by $\kappa$. It is clear that we can further identify $\cC$ with the quotient $\Az_f \times \Az_f / \sim$ with $(\ba,\bb) \sim (\bc,\bd)$ $\overset{def}{\LRarr}$ \an{$\ba \Rbar = \bc \Rbar$ and $\bb - \bd \in \ba \Rbar = \bc \Rbar$}. Let $\rho: \cC \to \Az_f \times \Az_f / \sim$, $\bb + \ba \Rbar \ma [\bb,\ba]$ denote this identification. Here $[\bb,\ba]$ is the equivalence class of $(\bb,\ba)$ in $\Az_f \times \Az_f / \sim$. Also, let $\pi: \Az_f \times \Az_f \onto \Az_f \times \Az_f / \sim$ be the canonical projection.
\blemma
\label{basic-open-sets}
The bijection $\rho \circ \kappa \circ \omega: \Spec(D_{\ping}) \to \Az_f \times \Az_f / \sim$ transports the topology of pointwise convergence on $\Spec(D_{\ping})$ to the quotient topology of $\Az_f \times \Az_f / \sim$ inherited from the natural topology of $\Az_f \times \Az_f$.
\elemma
\bproof
We have to prove that $\menge{(\rho \circ \kappa)(U(X; X_1, \dotsc, X_n))}{X, X_1, \dotsc, X_n \in \cJ_{\ping}}$ is a basis of open sets for the quotient topology on $\Az_f \times \Az_f / \sim$. We first compute
\bgloz
  \kappa(U(X; X_1, \dotsc, X_n)) 
  = \menge{\bb + \ba \Rbar}{\bb + \ba \Rbar \subseteq X(\Rbar), \bb + \ba \Rbar \nsubseteq X_i(\Rbar) \: \forall \: 1 \leq i \leq n}.
\egloz
Write $X(\Rbar) = b + \mfa \cdot \Rbar$, $X_i(\Rbar) = b_i + \mfa_i \cdot \Rbar$ for some $b$, $b_i$ in $K$ and fractional ideals $\mfa$, $\mfa_i$ of $K$ ($1 \leq i \leq n$). With this notation, $\rho(\kappa(U(X; X_1, \dotsc, X_n)))$ is the set of all $[\ba,\bb] \in \Az_f \times \Az_f / \sim$ satisfying
\begin{itemize}
\item $v_\mfp(\ba) \geq v_\mfp(\mfa)$, $v_\mfp(-b+\bb) \geq v_\mfp(\mfa)$ for all prime ideals $\mfp$;
\item for all $1 \leq i \leq n$, there exists a prime ideal $\mfp$ such that $v_\mfp(\ba) < v_\mfp(\mfa_i)$ or $v_\mfp(\ba) \geq v_\mfp(\mfa_i)$ and $v_\mfp(- b_i + \bb) < v_\mfp(\mfa_i)$.
\end{itemize}
From this description, we see that $\pi^{-1}(\rho(\kappa(U(X; X_1, \dotsc, X_n))))$ is clearly an open subset of $\Az_f \times \Az_f$. Moreover, every open subset of $\Az_f \times \Az_f$ of the form $\pi^{-1}(U)$ for some $U \subseteq \Az_f \times \Az_f / \sim$ is a union of sets of the form
\bgloz
  \rukl{\prod_{\mfp \notin \ti{F}} \mfp^{n_\mfp} \times \prod_{\mfp \in \ti{F}} (b_\mfp + \mfp^{n_\mfp})}
  \times
  \rukl{\prod_{\mfp \notin F} \mfp^{n_\mfp} \times \prod_{\mfp \in F} (\mfp^{n_\mfp} \setminus \mfp^{m_\mfp})} 
\egloz
for finite sets of prime ideals $F$ and $\ti{F}$, $b_\mfp \in K$ and integers $m_\mfp$, $n_\mfp$ which are $0$ for almost all prime ideals $\mfp$. As every such subset equals $\pi^{-1}(\rho(\kappa(U(X; X_1, \dotsc, X_n))))$ for some $X$, $X_1$, ..., $X_n$ in $\cJ_{\ping}$, we are done.
\eproof

\bcor
$\rho \circ \kappa \circ \omega: \Spec(D_{\ping}) \cong \Az_f \times \Az_f / \sim$ is a $G$-equivariant homeomorphism, where $G = K \rtimes K\reg$ acts on $\Az_f \times \Az_f / \sim$ by left multiplication, $(b,a) \cdot [\bb,\ba] = [b+a\bb,a\ba]$.
\ecor
\bcor
\label{A_f-D}
The transpose of $\rho \circ \kappa \circ \omega$ yields an identification $C_0(\Az_f \times \Az_f / \sim) \rtimes_r (K \rtimes K\reg) \cong D_{\ping} \rtimes_r G$.
\ecor

\section{Explicit description of primitive ideals}

For the sake of brevity, we write $P = R \rtimes R\reg$ and $G = K \rtimes K\reg$. The primitive ideal space of $C^*_r(P)$ has been computed in \cite{Ech-La}. Here is the result:

For every $\ba = (a_\mfp)_\mfp \in \Az_f$, set $Z(\ba) = \menge{\mfp \text{ prime ideal}}{a_\mfp = 0}$. Moreover, for a subset $A$ of the set of prime ideals $\cP$, let $V_A$ be the vanishing ideal of $C_0(\Az_f \times \Az_f / \sim)$ corresponding to the closed subset $C_A = \menge{[\bb,\ba] \in \Az_f \times \Az_f / \sim}{A \subseteq Z(\ba)}$, i.e. $V_A = \menge{f \in C_0(\Az_f \times \Az_f / \sim)}{f \vert_{C_A} = 0}$.
\btheo[\cite{Ech-La}]
\label{prim}
The map $2^{\cP} \ni A \ma \spkl{V_A} \in \Prim(C_0(\Az_f \times \Az_f / \sim) \rtimes_r (K \rtimes K\reg))$ is a bijection.
\etheo
\bproof
This follows from \cite[Corollary~2.7, Corollary~3.5]{Ech-La} using the easy observation that $\spkl{V_A} = V_A \rtimes_r (K \rtimes K\reg) = \ind V_A$.
\eproof
\bremark
As it is shown in \cite{Ech-La}, the map $\Az_f \times \Az_f / \sim \ \lori 2^{\cP}$, $[\bb,\ba] \ma Z(\ba)$ transports the quotient topology of the quasi-orbit space inherited from the topology of $\Az_f \times \Az_f / \sim$ to the power-cofinite topology of $2^{\cP}$. Basic open sets in this topology are of the form $\menge{\cS \in 2^{\cP}}{\cS \cap F = \emptyset}$, for $F \subseteq \cP$ finite. So with this topology on $2^{\cP}$, the bijection of the last corollary becomes a homeomorphism.
\eremark
\bcor
The ideals $\spkl{V_{\gekl{\mfp}}}$, $\mfp \in \cP$, are the minimal non-zero primitive ideals of $C_0(\Az_f \times \Az_f / \sim) \rtimes_r (K \rtimes K\reg)$.
\ecor
\bproof
The bijection $2^{\cP} \ni A \ma \spkl{V_A} \in \Prim(C_0(\Az_f \times \Az_f / \sim) \rtimes_r (K \rtimes K\reg))$ is clearly inclusion-preserving.
\eproof

Theorem~\ref{prim} gives a description of the primitive ideal space of $C_0(\Az_f \times \Az_f / \sim) \rtimes_r (K \rtimes K\reg)$, hence also of $D_{\ping} \rtimes_r G$. In the sequel, we describe the minimal non-zero primitive ideals of $D_{\ping} \rtimes_r G$ in an explicit way. Such a description could also be obtained for all the primitive ideals, but it is really the minimal ones we will be interested in later on.

Let $J_{\gekl{\mfp}}$ be the ideal of $D_{\ping} \rtimes_r G$ which corresponds to the ideal $\spkl{V_{\gekl{\mfp}}}$ of $C_0(\Az_f \times \Az_f / \sim) \rtimes_r (K \rtimes K\reg)$ under the isomorphism $C_0(\Az_f \times \Az_f / \sim) \rtimes_r (K \rtimes K\reg) \cong D_{\ping} \rtimes_r G$ from Corollary~\ref{A_f-D}.
\bprop
$J_{\gekl{\mfp}} = \spkl{E_{R \times R\reg} - E_{R \times \mfp\reg}}_{D_{\ping} \rtimes_r G}$.
\eprop
\bproof
We first show \an{$\supseteq$}. Take $[\bb,\ba] \in C_{\gekl{\mfp}}$, i.e. $a_\mfp = 0$ ($\ba=(a_\mfp)_{\mfp}$), and set $\chi = (\rho \circ \kappa \circ \omega)^{-1}[\bb,\ba]$. If $\bb + \ba \Rbar \nsubseteq \Rbar$, then for all $x \in R$, $\bb + \ba \Rbar \nsubseteq x + \mfp \cdot \Rbar$, hence $\chi(E_{(x+\mfp) \times \mfp \reg}) = 0$. As $R \times \mfp\reg = \bigcup_{x \in R} (x+\mfp) \times \mfp\reg$, we deduce $\chi(E_{R \times R\reg}) = \chi(E_{R \times \mfp\reg}) = 0$. If $\bb + \ba \Rbar \subseteq \Rbar$, then since $a_\mfp = 0$, there exists a fractional ideal $\mfa$ such that $v_\mfp(\mfa) \geq 1$ and $\bb + \ba \Rbar \subseteq b + \mfa \cdot \Rbar$ for some $b \in K$. Then $\bb + \ba \Rbar \subseteq \Rbar \cap (b + \mfa \cdot \Rbar) = x + (\mfa \cap R) \cdot \Rbar$ for some $x \in R$. As $v_\mfp(\mfa) \geq 1$, we deduce $\mfa \cap R \subseteq \mfp$ and thus $\bb + \ba \Rbar \subseteq x + \mfp \cdot \Rbar \subseteq \Rbar$. This implies $\chi(E_{R \times R\reg}) = \chi(E_{(x+\mfp) \times \mfp\reg}) = 1$. Since $\menge{E_{(r+\mfp) \times \mfp\reg}}{r \in R}$ are pairwise orthogonal and $E_{R \times \mfp\reg} = \sum_{r \in R} E_{(r+\mfp) \times \mfp\reg}$, we conclude that $\chi(E_{R \times \mfp\reg}) = 1$. Hence $\chi(E_{R \times R\reg} - E_{R \times \mfp\reg}) = 0$.

To prove the reverse inclusion, take $[\bb,\ba] \in \Az_f \times \Az_f / \sim$ and set $\chi = (\rho \circ \kappa \circ \omega)^{-1}[\bb,\ba]$. We need to show that if $\chi(\spkl{E_{R \times R\reg} - E_{R \times \mfp\reg}}_{D_{\ping} \rtimes_r G} \cap D_{\ping}) = 0$, then $a_\mfp = 0$, where $\ba = (a_{\mfp})_{\mfp}$. $a_\mfp = 0$ means that
\bgloz
  \sup \menge{v_\mfp(\mfa)}{\text{There is } (b+\mfa) \times \mfa\reg \in \cJ_{\ping} \text{ with } \chi(E_{(b+\mfa) \times \mfa\reg}) = 1} = \infty.
\egloz
Assume that this supremum is finite. Then we can choose $b \in K$ and a fractional ideal $\mfa$ of $K$ such that $\chi(E_{(b+\mfa) \times \mfa\reg}) = 1$ and the supremum above agrees with $v_\mfp(\mfa)$. By strong approximation, we can find $a \in K\reg$ such that $v_\mfp(a) = v_\mfp(\mfa)$ and $\mfa \subseteq aR$. The latter condition implies $\chi(E_{(b+aR) \times (aR)\reg}) = 1$. Therefore, $\chi(\Ad(U_{(b,a)})(E_{R \times R\reg})) = 1$. As $\chi$ vanishes on $\spkl{E_{R \times R\reg} - E_{R \times \mfp\reg}}_{D_{\ping} \rtimes_r G} \cap D_{\ping}$, it follows that
\bgloz
  0 = \chi \circ \Ad(U_{(b,a)}) (E_{R \times R\reg} - E_{R \times \mfp\reg}) = 1 - \sum_{x \in R} \chi \circ \Ad(U_{(b,a)}) (E_{(x+\mfp) \times \mfp\reg}).
\egloz
Thus we must have $1 = \chi \circ \Ad(U_{(b,a)}) (E_{(x+\mfp) \times \mfp\reg}) = \chi(E_{(b+ax+a\mfp) \times (a\mfp)\reg})$ for some $x \in R$. But $v_\mfp(a\mfp) = v_\mfp(a) + 1 = v_\mfp(\mfa) + 1 > v_\mfp(\mfa)$. This contradicts maximality of $v_\mfp(\mfa)$.
\eproof

In the sequel, we describe $J_{\gekl{\mfp}}$ in a way which is more suitable for K-theoretic computations.

\blemma
The smallest family $\cJ(\gekl{\mfp})$ of subsets of $G = K \rtimes K\reg$ with the properties
\begin{itemize}
\item $R \times (R \setminus \mfp) \in \cJ(\gekl{\mfp})$,
\item $g \in G, X \in \cJ(\gekl{\mfp}) \Rarr g \cdot X \in \cJ(\gekl{\mfp})$,
\item $X_1, X_2 \in \cJ(\gekl{\mfp}) \Rarr X_1 \cap X_2 \in \cJ(\gekl{\mfp})$,
\item $X \in \cJ(\gekl{\mfp}), Y \in \cJ_{\ping} \Rarr X \cap Y \in \cJ(\gekl{\mfp})$
\end{itemize}
is given by
\bgl
\label{set}
  \menge{(b+\mfa) \times (\mfa \setminus \mfp \cdot \mfa)}{b \in K, \mfa \text{ a fractional ideal of } K} \cup \gekl{\emptyset}.
\egl
\elemma
\bproof
The first item is obviously fulfilled. And the set in \eqref{set} is certainly $G$-invariant, i.e. satisfies the second item. To prove the third property, take two fractional ideals $\mfa = \mfp^m \cdot \mfp_1^{m_1} \dotsm \mfp_r^{m_r}$, $\mfb = \mfp^n \cdot \mfp_1^{n_1} \dotsm \mfp_r^{n_r}$ with prime ideals $\mfp_i$. Then $(\mfa \setminus \mfp \cdot \mfa) \cap (\mfb \setminus \mfp \cdot \mfb) = ((\mfp^m \setminus \mfp^{m+1}) \cap (\mfp^n \setminus \mfp^{n+1})) \cdot \mfp_1^{l_1} \dotsm \mfp_r^{l_r}$ with $l_i \defeq \max(m_i,n_i)$. As $(\mfp^m \setminus \mfp^{m+1}) \cap (\mfp^n \setminus \mfp^{n+1})$ is either empty or $\mfp^m \setminus \mfp^{m+1}$ depending on whether $m \neq n$ or $m=n$, we see that $(\mfa \setminus \mfp \cdot \mfa) \cap (\mfb \setminus \mfp \cdot \mfb)$ is either empty or $\mfp^m \cdot \mfp_1^{l_1} \dotsm \mfp_r^{l_r} \setminus \mfp^{m+1} \cdot \mfp_1^{l_1} \dotsm \mfp_r^{l_r} = (\mfa \cap \mfb) \setminus (\mfp \cdot (\mfa \cap \mfb))$. For the fourth property, let $\mfa$ and $\mfb$ be as above. Then $(\mfa \setminus \mfp \cdot \mfa) \cap \mfb = ((\mfp^m \setminus \mfp^{m+1}) \cap \mfp^n) \cdot \mfp_1^{l_1} \dotsm \mfp_r^{l_r}$ with $l_i \defeq \max(m_i,n_i)$. Thus either $(\mfa \setminus \mfp \cdot \mfa) \cap \mfb = \emptyset$ or $(\mfa \setminus \mfp \cdot \mfa) \cap \mfb = (\mfp^m \setminus \mfp^{m+1}) \cdot \mfp_1^{l_1} \dotsm \mfp_r^{l_r} = \mfp^m \cdot \mfp_1^{l_1} \dotsm \mfp_r^{l_r} \setminus \mfp^{m+1} \cdot \mfp_1^{l_1} \dotsm \mfp_r^{l_r}$ depending on whether $m<n$ or $m \geq n$. So far, we have proven that the set in \eqref{set} satisfies all the desired properties. It remains to prove minimality. It suffices to show that for every non-zero ideal $\mfa$ of $R$, we must have $\mfa \times (\mfa \setminus \mfp \cdot \mfa) \in \cJ(\gekl{\mfp})$. Write $\mfa = \mfp^v \cdot \mfa'$ for some $v \in \Nz_0$ and some non-zero ideal $\mfa'$ of $R$ which is coprime to $\mfp$. As every fractional ideal is of the form $(aR) \cap (cR)$ for some $a$, $c$ in $K\reg$, we can choose $a, c \in K\reg$ with $\mfp^v = (aR) \cap (cR)$, and we can without loss of generality assume $v_\mfp(a) = v$. Then $a(R \setminus \mfp) \cap cR = (aR \cap cR) \setminus (a\mfp \cap cR) = \mfp^v \setminus \mfp^{v+1}$. Therefore
\bglnoz
  && \mfa \times (\mfa \setminus \mfp \cdot \mfa) = (\mfa' \times (\mfa')\reg) \cap (\mfp^v \times (\mfp^v \setminus \mfp^{v+1})) \\
  &=& (\mfa' \times (\mfa')\reg) \cap (aR \times (a(R \setminus \mfp)) \cap (cR \times (cR)\reg) \in \cJ(\gekl{\mfp}).
\eglnoz
\eproof

For a fractional ideal $\mfa$, set $\delta_{\mfa,\mfp} = E_{\mfa \times \mfa \reg} - E_{\mfa \times (\mfp \cdot \mfa)\reg} = E_{\mfa \times (\mfa \setminus \mfp \cdot \mfa)}$. With this notation, we have $J_{\gekl{\mfp}} = \spkl{\delta_{R,\mfp}}$.
\bcor
The smallest $G$-invariant ideal of $D_{\ping}$ containing $\delta_{R,\mfp}$ is given by $C^*(\menge{E_X}{X \in \cJ(\gekl{\mfp})})$.
\ecor
\bcor
We have $\spkl{\delta_{R,\mfp}}_{D_{\ping} \rtimes_r G} = C^*(\menge{E_X}{X \in \cJ(\gekl{\mfp})}) \rtimes_r G$.
\ecor

For our K-theoretic computations, we need
\blemma
\label{J-ind}
$\cJ(\gekl{\mfp})$ is independent.
\elemma
\bproof
Assume that $(b+\mfa) \times (\mfa \setminus \mfp \cdot \mfa) = \bigcup_{i=1}^n (b_i + \mfa_i) \times (\mfa_i \setminus \mfp \cdot \mfa_i)$. Thus $\mfa \setminus \mfp \cdot \mfa = \bigcup_{i=1}^n (\mfa_i \setminus \mfp \cdot \mfa_i)$. Moreover, we deduce that $\mfa\reg \supseteq \mfa \setminus \mfp \cdot \mfa \supseteq \mfa_i \setminus \mfp \cdot \mfa_i$. This implies $\mfa\reg \cup (\mfp \cdot \mfa_i)\reg \supseteq \mfa_i\reg$. Using strong approximation as in \cite[Lemma~2.30]{Li1}, we deduce that $\mfa_i \subseteq \mfa$ ($*$). Hence, $\mfa\reg = \bigcup_{i=1}^n (\mfa_i \setminus \mfp \cdot \mfa_i) \cup (\mfp \cdot \mfa)\reg \overset{(*)}{=} \bigcup_{i=1}^n \mfa_i \reg \cup (\mfp \cdot \mfa)\reg$. Again, as in \cite[Lemma~2.30]{Li1}, strong approximation implies that there exists $1 \leq i \leq n$ with $\mfa = \mfa_i$.
\eproof

All in all, we obtain
\bprop
\label{minprimid}
For every prime ideal $\mfp$, let $\delta_{R,\mfp} = E_{R \times R\reg} - E_{R \times \mfp\reg}$ (with $\mfp\reg = \mfp \setminus \gekl{0}$). Moreover, set 
\bgloz
  \cJ(\gekl{\mfp}) \defeq \menge{(b+\mfa) \times (\mfa \setminus \mfp \cdot \mfa)}{b \in K, \mfa \text{ a fractional ideal of } K} \cup \gekl{\emptyset}.
\egloz
Then the minimal non-zero primitive ideals of $D_{\ping} \rtimes_r G$ are given by
\bgloz
  J_{\gekl{\mfp}} = \spkl{\delta_{R,\mfp}} = C^*(\menge{E_X}{X \in \cJ(\gekl{\mfp})}) \rtimes_r G, \ \mfp \in \cP.
\egloz
Moreover, $\cJ(\gekl{\mfp})$ is independent.
\eprop

\section{K-theoretic invariants}
\label{K-thinv}

Our goal in this section is to compute the K-theory of minimal non-zero primitive ideals of $C^*_r(R \rtimes R\reg)$. Moreover, we also determine the torsion order of the class of the unit in the corresponding quotients. As before, we write $P = R \rtimes R\reg$ and $G = K \rtimes K\reg$.

The K-theory of $C^*_r(R \rtimes R\reg)$ has already been computed in \cite{C-E-L1}. Let us recall the result. For every class $\mfk$ in the class group $Cl_K$, let $\mfa_{\mfk}$ be an ideal of $R$ representing $\mfk$. If $\mfk$ is the trivial class in $Cl_K$, we choose $\mfa_{\mfk} = R$.

\btheooz[Theorem~8.2.1 in \cite{C-E-L1}]
The homomorphisms $\iota_{\mfk}: C^*(\mfa_{\mfk} \rtimes R^*) \to C^*_r(R \rtimes R\reg)$ given by $u_x \ma V_x E_{\mfa_{\mfk} \times \mfa_{\mfk}\reg}$ ($\mfa_{\mfk}\reg = \mfa_{\mfk} \setminus \gekl{0}$) induce an isomorphism
\bgloz
  \sum_{\mfk} (\iota_{\mfk})_*: \bigoplus_{\mfk \in Cl_K} K_*(C^*(\mfa_{\mfk} \rtimes R^*)) \cong K_*(C^*_r(R \rtimes R\reg)).
\egloz
\etheooz
Here and in the sequel, the $u_x$ are the canonical unitaries in the group C*-algebra $C^*(\mfa_{\mfk} \rtimes R^*)$.

As above, we write $\cP$ for the set of non-zero prime ideals of $R$, and we set $\delta_{\mfa,\mfp} = E_{\mfa \times \mfa\reg} - E_{\mfa \times (\mfp \cdot \mfa)\reg}$ for every $(0) \neq \mfa \triangleleft R$ and $\mfp \in \cP$.
\btheo
The minimal non-zero primitive ideals of $C^*_r(R \rtimes R\reg)$ are given by $I_{\gekl{\mfp}} \defeq \spkl{\delta_{R,\mfp}}_{C^*_r(R \rtimes R\reg)}$, $\mfp \in \cP$. Moreover, for every prime ideal $\mfp$, the homomorphisms $\kappa_{\mfk}: C^*(\mfa_{\mfk} \rtimes R^*) \to I_{\gekl{\mfp}}$, $u_x \ma V_x \delta_{\mfa_{\mfk},\mfp}$ induce an isomorphism
\bgloz
  \sum_{\mfk} (\kappa_{\mfk})_*: \bigoplus_{\mfk \in Cl_K} K_*(C^*(\mfa_{\mfk} \rtimes R^*)) \cong K_*(I_{\gekl{\mfp}}).
\egloz
\etheo
\bproof
Since $C^*_r(R \rtimes R\reg)$ is isomorphic to a full corner of $D_{\ping} \rtimes_r G$ for $P = R \rtimes R\reg$, $G = K \rtimes K\reg$, we know that minimal non-zero primitive ideals of $C^*_r(R \rtimes R\reg)$ are in one-to-one correspondence with minimal non-zero primitive ideals of $D_{\ping} \rtimes_r G$. For every prime ideal $\mfp$, the ideal $J_{\gekl{\mfp}}$ of $D_{\ping} \rtimes_r G$ corresponds to the ideal $E_P J_{\gekl{\mfp}} E_P$ of $E_P (D_{\ping} \rtimes_r G) E_P$, hence to the ideal $I_{\gekl{\mfp}} = \spkl{\delta_{R,\mfp}}_{C^*_r(R \rtimes R\reg)}$ of $C^*_r(R \rtimes R\reg)$. This together with Proposition~\ref{minprimid} proves the first part of our theorem.

For the K-theoretic formula, just apply \cite[Corollary~3.14]{C-E-L2} with $G = K \rtimes K\reg$ and $I = \cJ(\gekl{\mfp})\reg = \cJ(\gekl{\mfp}) \setminus \gekl{\emptyset}$.
\eproof

Let us study the inclusion $i_{\gekl{\mfp}}: I_{\gekl{\mfp}} \into C^*_r(R \rtimes R\reg)$ in K-theory. Recall that $\iota_{\mfk}: C^*(\mfa_{\mfk} \rtimes R^*) \to C^*_r(R \rtimes R\reg)$, $u_x \ma V_x E_{\mfa_{\mfk} \times \mfa_{\mfk}\reg}$ and $\kappa_{\mfk}: C^*(\mfa_{\mfk} \rtimes R^*) \to I_{\gekl{\mfp}}$, $u_x \ma V_x \delta_{\mfa_{\mfk},\mfp}$ are the homomorphisms from the theorems above. Moreover, we introduce another homomorphism $\iota_{\mfk,\mfp}: C^*(\mfa_{\mfk} \rtimes R^*) \to C^*_r(R \rtimes R\reg)$, $u_x \ma V_x E_{\mfa_{\mfk} \times (\mfp \cdot \mfa_{\mfk})\reg}$. It is clear that in K-theory, $i_{\gekl{\mfp}} \circ \kappa_{\mfk}$ is the difference of the two homomorphisms $\iota_{\mfk}$ and $\iota_{\mfk,\mfp}$. To describe $\iota_{\mfk,\mfp}$ in $K_0$, we choose for every $\mfk \in Cl_K$ an element $a(\mfp \cdot \mfa_{\mfk}) \in K\reg$ such that $a(\mfp \cdot \mfa_{\mfk}) \cdot \mfp \cdot \mfa_{\mfk} = \mfa_{[\mfp] \cdot \mfk}$. This element $a(\mfp \cdot \mfa_{\mfk})$ induces the following isomomorphism:
\bgloz
  \mu_{a(\mfp \cdot \mfa_{\mfk})}: C^*((\mfp \cdot \mfa_{\mfk}) \rtimes R^*) \cong C^*(\mfa_{[\mfp] \cdot \mfk} \rtimes R^*), \ 
  u_{(b,a)} \ma u_{(a(\mfp \cdot \mfa_{\mfk}) b,a)}.
\egloz
The element $a(\mfp \cdot \mfa_{\mfk})$ is determined up to a unit in $R\reg$, so that $\mu_{a(\mfp \cdot \mfa_{\mfk})}$ is determined up to an inner automorphism. Hence in K-theory, $(\mu_{a(\mfp \cdot \mfa_{\mfk})})_*$ does not depend on the choice of $a(\mfp \cdot \mfa_{\mfk})$.

\blemma
$(\iota_{\mfk,\mfp})_* 
= (\iota_{[\mfp] \cdot \mfk})_* \circ (\mu_{a(\mfp \cdot \mfa_{\mfk})})_* \circ \res_{\mfa_{\mfk} \rtimes R^*}^{(\mfp \cdot \mfa_{\mfk}) \rtimes R^*}$.
\elemma
\bproof
The image of $\iota_{\mfk,\mfp}$ lies in
\bgloz
  C_{\mfk,\mfp} \defeq 
  C^*(\menge{E_{(r + \mfp \cdot \mfa_{\mfk}) \times (\mfp \cdot \mfa_{\mfk})\reg}}{r \in \mfa_{\mfk} / \mfp \cdot \mfa_{\mfk}} 
  \cup 
  \menge{V_x E_{\mfa_{\mfk} \times (\mfp \cdot \mfa_{\mfk}) \reg}}{x \in \mfa_{\mfk} \rtimes R^*}).
\egloz
It is clear that $C_{\mfk,\mfp}$ can be identified with $C(\mfa_{\mfk} / \mfp \cdot \mfa_{\mfk}) \rtimes (\mfa_{\mfk} \rtimes R^*) \cong C(\mfa_{\mfk} \rtimes R^* / (\mfp \cdot \mfa_{\mfk}) \rtimes R^*) \rtimes (\mfa_{\mfk} \rtimes R^*)$ in such a way that $\iota_{\mfk,\mfp}$ corresponds to the canonical embedding $C^*(\mfa_{\mfk} \rtimes R^*) \to C(\mfa_{\mfk} \rtimes R^* / (\mfp \cdot \mfa_{\mfk}) \rtimes R^*) \rtimes (\mfa_{\mfk} \rtimes R^*)$. Now let $\rho_{\mfk,\mfp}$ denote the composition $C_{\mfk,\mfp} \cong C(\mfa_{\mfk} \rtimes R^* / (\mfp \cdot \mfa_{\mfk}) \rtimes R^*) \rtimes (\mfa_{\mfk} \rtimes R^*) \cong M_{N(\mfp)}(C^*((\mfp \cdot \mfa_{\mfk}) \rtimes R^*))$. The last isomorphism is the homomorphism $\phi$ from the appendix (defined before \eqref{phi-res}), for $G = \mfa_{\mfk} \rtimes R^*$ and $H = (\mfp \cdot \mfa_{\mfk}) \rtimes R^*$. Now consider the diagram
\bgl
\label{CD}
  \xymatrix@C=1mm{
  & M_{N(\mfp)}(C^*((\mfp \cdot \mfa_{\mfk}) \rtimes R^*)) & \ \ \ \ \ \ & 
  \ar[ll]_{\ \ \ e \otimes \id} C^*((\mfp \cdot \mfa_{\mfk}) \rtimes R^*) 
  \ar[d]_{\cong}^{\mu_{a(\mfp \cdot \mfa_{\mfk})}} \\
  C^*(\mfa_{\mfk} \rtimes R^*) 
  \ar[r]^{\iota_{\mfk,\mfp} \vert^{C_{\mfk,\mfp}}} \ar[dr]^{\iota_{\mfk,\mfp}} & C_{\mfk,\mfp} \ar[u]^{\cong}_{\rho_{\mfk,\mfp}} \ar[d]^{i_{C_{\mfk,\mfp}}} 
  & & \ar[dll]^{\iota_{[\mfp] \cdot \mfk}} C^*(\mfa_{[\mfp] \cdot \mfk} \rtimes R^*) \\
  & C^*_r(R \rtimes R\reg) & &
  }
\egl
where $i_{C_{\mfk,\mfp}}$ is the canonical inclusion and $e \otimes \id$ is the canonical embedding into the upper left corner. The lower left triangle in \eqref{CD} commutes by construction. The pentagon on the right of \eqref{CD} commutes in K-theory because once we compose with the inclusion $C^*_r(R \rtimes R\reg) \into D_{\ping} \rtimes_r G$, the homomorphisms $\iota_{[\mfp] \cdot \mfk} \circ \mu_{a(\mfp \cdot \mfa_{\mfk})}$ and $i_{C_{\mfk,\mfp}} \circ (\rho_{\mfk,\mfp})^{-1} \circ (e \otimes \id)$ only differ by an inner automorphism. Therefore, we obtain
\bgloz
  (\iota_{\mfk,\mfp})_* = (i_{C_{\mfk,\mfp}})_* \circ (\iota_{\mfk,\mfp} \vert^{C_{\mfk,\mfp}})_*
  = (\iota_{[\mfp] \cdot \mfk})_* \circ (\mu_{a(\mfp \cdot \mfa_{\mfk})})_* \circ (e \otimes \id)_*^{-1} \circ (\rho_{\mfk,\mfp})_* 
  \circ (\iota_{\mfk,\mfp} \vert^{C_{\mfk,\mfp}})_*.
\egloz
But $\rho_{\mfk,\mfp} \circ \iota_{\mfk,\mfp} \vert^{C_{\mfk,\mfp}}$ coincides with the homomorphism $\varphi$ from \eqref{phi-res} in the appendix for $G = \mfa_{\mfk} \rtimes R^*$, $H = (\mfp \cdot \mfa_{\mfk}) \rtimes R^*$. Thus by Lemma~\ref{lem-res1} in the appendix, $(e \otimes \id)_*^{-1} \circ (\rho_{\mfk,\mfp})_* \circ (\iota_{\mfk,\mfp} \vert^{C_{\mfk,\mfp}})_* = \res_{\mfa_{\mfk} \rtimes R^*}^{(\mfp \cdot \mfa_{\mfk}) \rtimes R^*}$ and we are done.
\eproof

For every prime ideal $\mfp$, let $h_{\mfp}$ be the order of $[\mfp]$ in $Cl_K$ and let $a(\mfp^{h_{\mfp}}) \in K\reg$ satisfy $a(\mfp^{h_{\mfp}}) \cdot \mfp^{h_{\mfp}} = R$.
\blemma
\label{imim}
Upon identifying $K_0(C^*(R \rtimes R^*))$ with the corresponding direct summand in $\bigoplus_{\mfk \in Cl_K} K_0(C^*(\mfa_{\mfk} \rtimes R^*))$, we have that
\bgloz
  \img \rukl{(\sum_{\ti{\mfk}} (\iota_{\ti{\mfk}})_*)^{-1} \circ i_{\gekl{\mfp}} \circ (\sum_{\mfk} (\kappa_{\mfk})_*)} \cap \Zz [1]_{[R]}
\egloz
and
\bgloz
  (\id - (\mu_{a(\mfp^{h_{\mfp}})})_* \circ \res_{R \rtimes R^*}^{\mfp^{h_{\mfp}} \rtimes R^*})(K_0(C^*(R \rtimes R^*))) \cap \Zz [1]
\egloz
coincide. Here $[1]$ is the class of the unit of $C^*(R \rtimes R^*)$ and $[1]_{[R]}$ is the image of $[1]$ in $\bigoplus_{\mfk \in Cl_K} K_0(C^*(\mfa_{\mfk} \rtimes R^*))$ under the canonical inclusion $K_0(C^*(R \rtimes R^*)) \into \bigoplus_{\mfk} K_0(C^*(\mfa_{\mfk} \rtimes R^*))$.
\elemma
\bproof
The previous lemma tells us that
\bgl
\label{ik=i-i}
  (i_{\gekl{\mfp}} \circ \kappa_{\gekl{\mfp}})_* = (\iota_{\mfk})_* - (\iota_{\mfk,\mfp})_* 
  = (\iota_{\mfk})_* - (\iota_{[\mfp] \cdot \mfk})_*
  \circ (\mu_{a(\mfp \cdot \mfa_{\mfk})})_* \circ \res_{\mfa_{\mfk} \rtimes R^*}^{(\mfp \cdot \mfa_{\mfk}) \rtimes R^*}.
\egl
Therefore, it suffices to consider only those $\mfk \in Cl_K$ which lie in the subgroup of $Cl_K$ generated by $[\mfp]$, i.e. $\spkl{[\mfp]} = \gekl{[R],[\mfp], \dotsc, [\mfp^{h_{\mfp}}]}$. In other words,
\bglnoz
  && \img \rukl{(\sum_{\ti{\mfk}} (\iota_{\ti{\mfk}})_*)^{-1} \circ i_{\gekl{\mfp}} \circ (\sum_{\mfk} (\kappa_{\mfk})_*)} \cap \Zz [1]_{[R]} \\
  &=& \img \rukl{(\sum_{\ti{h} = 0}^{h_{\mfp}-1} (\iota_{[\mfp^{\ti{h}}]})_*)^{-1} \circ i_{\gekl{\mfp}} \circ (\sum_{h=0}^{h_{\mfp}-1} (\kappa_{[\mfp^h]})_*)} 
  \cap \Zz [1]_{[R]}.
\eglnoz
Using \eqref{ik=i-i}, we further compute
\bglnoz
  && (\sum_{\ti{h} = 0}^{h_{\mfp}-1} (\iota_{[\mfp^{\ti{h}}]})_*)^{-1} \circ i_{\gekl{\mfp}} \circ (\sum_{h=0}^{h_{\mfp}-1} (\kappa_{[\mfp^h]})_*) \\
  &=& (\sum_{\ti{h}} (\iota_{[\mfp^{\ti{h}}]})_*)^{-1} \circ \rukl{\sum_h (\iota_{[\mfp^{h}]})_* 
  - (\iota_{[\mfp^{h+1}]})_* \circ (\mu_{a(\mfp \cdot \mfa_{[\mfp^h]})})_* \circ \res_{\mfa_{[\mfp^h] \rtimes R^*}}^{(\mfp \cdot \mfa_{[\mfp^h]}) \rtimes R^*}}.
\eglnoz
Let $\ve_{[\mfp^h]}$ be the canonical embedding $K_0(C^*(\mfa_{[\mfp^h]} \rtimes R^*)) \to \bigoplus_{\mfk \in \spkl{[\mfp]}} K_0(C^*(\mfa_{\mfk} \rtimes R^*))$. Then
\bgloz
  (\sum_{\ti{h} = 0}^{h_{\mfp}-1} (\iota_{[\mfp^{\ti{h}}]})_*)^{-1} \circ i_{\gekl{\mfp}} \circ (\sum_{h=0}^{h_{\mfp}-1} (\kappa_{[\mfp^h]})_*)
  = \sum_h \rukl{\ve_{[\mfp^{h}]} - \ve_{[\mfp^{h+1}]} \circ (\mu_{a(\mfp \cdot \mfa_{[\mfp^h]})})_* \circ 
  \res_{\mfa_{[\mfp^h]} \rtimes R^*}^{(\mfp \cdot \mfa_{[\mfp^h]}) \rtimes R^*}}
\egloz
on $\bigoplus_{\mfk \in \spkl{[\mfp]}} K_0(C^*(\mfa_{\mfk} \rtimes R^*)) \subseteq \bigoplus_{\mfk \in Cl_K} K_0(C^*(\mfa_{\mfk} \rtimes R^*))$.

Now let $x_h$, $0 \leq h \leq h_{\mfp} - 1$, be elements of $K_0(C^*(\mfa_{[\mfp^h]} \rtimes R^*))$ such that
\bgloz
  \sum_h \rukl{\ve_{[\mfp^{h}]} - 
  \ve_{[\mfp^{h+1}]} \circ (\mu_{a(\mfp \cdot \mfa_{[\mfp^h]})})_* \circ \res_{\mfa_{[\mfp^h]} \rtimes R^*}^{(\mfp \cdot \mfa_{[\mfp^h]}) \rtimes R^*}}(x_h)
  = z [1]_{[R]} \in \Zz [1]_{[R]}.
\egloz
Since
\bglnoz
  && \sum_h \rukl{\ve_{[\mfp^{h}]} - 
  \ve_{[\mfp^{h+1}]} \circ (\mu_{a(\mfp \cdot \mfa_{[\mfp^h]})})_* \circ \res_{\mfa_{[\mfp^h]} \rtimes R^*}^{(\mfp \cdot \mfa_{[\mfp^h]}) \rtimes R^*}}(x_h) \\
  &=& \ve_{[R]} 
  \rukl{x_0 - 
  ((\mu_{a(\mfp \cdot \mfa_{[\mfp^{h_{\mfp}-1}]})})_* 
  \circ \res_{\mfa_{[\mfp^{h_{\mfp}-1}]} \rtimes R^*}^{(\mfp \cdot \mfa_{[\mfp^{h_{\mfp}-1}]} \rtimes R^*)})(x_{h_{\mfp}-1})}
  \\
  &+& \sum_{h=0}^{h_{\mfp}-2} \ve_{[\mfp^{h+1}]} \rukl{x_{h+1} 
  - ((\mu_{a(\mfp \cdot \mfa_{[\mfp^h]})})_* \circ \res_{\mfa_{[\mfp^h]} \rtimes R^*}^{(\mfp \cdot \mfa_{[\mfp^h]}) \rtimes R^*})(x_h)},
\eglnoz
we deduce that for every $0 \leq h \leq h_{\mfp} - 2$, we must have
\bgloz
  x_{h+1} = ((\mu_{a(\mfp \cdot \mfa_{[\mfp^h]})})_* \circ \res_{\mfa_{[\mfp^h]} \rtimes R^*}^{(\mfp \cdot \mfa_{[\mfp^h]}) \rtimes R^*})(x_h).
\egloz
Thus
\bgl
\label{z1}
  z [1]_{[R]} = x_0 - (\mu_{a(\mfp \cdot \mfa_{[\mfp^{h_{\mfp}-1}]})})_* 
  \circ \res_{\mfa_{[\mfp^{h_{\mfp}-1}]} \rtimes R^*}^{(\mfp \cdot \mfa_{[\mfp^{h_{\mfp}-1}]}) \rtimes R^*} 
  \circ \dotsb \circ (\mu_{a(\mfp)})_* \circ \res_{R \rtimes R^*}^{\mfp \rtimes R^*} (x_0).
\egl
As
\bgloz
  a(\mfp \cdot \mfa_{[\mfp^{h_{\mfp}-1}]}) \cdot a(\mfp \cdot \mfa_{[\mfp^{h_{\mfp}-2}]}) \cdots a(\mfp \cdot \mfa_{[\mfp]}) a(\mfp) \mfp^{h_{\mfp}} = R,
\egloz
we deduce that $a(\mfp \cdot \mfa_{[\mfp^{h_{\mfp}-1}]}) \cdot a(\mfp \cdot \mfa_{[\mfp^{h_{\mfp}-2}]}) \cdots a(\mfp \cdot \mfa_{[\mfp]}) a(\mfp)$ and $a(\mfp^{h_{\mfp}})$ only differ by a unit in $R\reg$. Thus
\bgloz
  (\mu_{a(\mfp \cdot \mfa_{[\mfp^{h_{\mfp}-1}]})})_* \circ \cdots (\mu_{a(\mfp)})_* = (\mu_{a(\mfp^{h_{\mfp}})})_*.
\egloz
By Lemma~\ref{lem-res2} in the appendix, this tells us that
\bgloz
  (\mu_{a(\mfp \cdot \mfa_{[\mfp^{h_{\mfp}-1}]})})_* \circ \res_{\mfa_{[\mfp^{h_{\mfp}-1}]} \rtimes R^*}^{(\mfp \cdot \mfa_{[\mfp^{h_{\mfp}-1}]}) \rtimes R^*}
  \circ \dotsb \circ (\mu_{a(\mfp)})_* \circ \res_{R \rtimes R^*}^{\mfp \rtimes R^*} 
  = (\mu_{a(\mfp^{h_{\mfp}})})_* \circ \res_{R \rtimes R^*}^{\mfp^{h_{\mfp}} \rtimes R^*}.
\egloz
Thus $z [1]_{[R]} = (\id - (\mu_{a(\mfp^{h_{\mfp}})})_* \circ \res_{R \rtimes R^*}^{\mfp^{h_{\mfp}} \rtimes R^*})(x_0)$.
\eproof

Now write $R^* = \mu \times \Gamma$ where $\mu$ is the group of roots of unity in $K$, and $\Gamma$ is a free abelian group. Let $1$ be the unit of $C^*(R \rtimes \Gamma)$. Fix an element $a(\mfp^{h_{\mfp}}) \in K\reg$ such that $a(\mfp^{h_{\mfp}}) \mfp^{h_{\mfp}} = R$. Moreover, let $\mu^{\Gamma}_{a(\mfp^{h_{\mfp}})}$ be the homomorphism $C^*(\mfp^{h_{\mfp}} \rtimes \Gamma) \cong C^*(R \rtimes \Gamma)$, $u_{(b,a)} \ma u_{(a(\mfp^{h_{\mfp}}) b,a)}$.
\blemma
$(\id - (\mu^{\Gamma}_{a(\mfp^{h_{\mfp}})})_* \circ \res_{R \rtimes \Gamma}^{\mfp^{h_{\mfp}} \rtimes \Gamma})(K_0(C^*(R \rtimes \Gamma))) \cap \Zz[1] = ((N(\mfp)^{h_{\mfp}} - 1)\Zz) [1]$.
\elemma
\bproof
Let $\tau^{R \rtimes \Gamma}$ and $\tau^{\mfp^{h_{\mfp}} \rtimes \Gamma}$ be the canonical tracial states on $C^*(R \rtimes \Gamma)$ and $C^*(\mfp^{h_{\mfp}} \rtimes \Gamma)$. We write $\tau^{R \rtimes \Gamma}_*$ and $\tau^{\mfp^{h_{\mfp}} \rtimes \Gamma}_*$ for the induced homomorphisms on $K_0$. As $R \rtimes \Gamma$ is amenable, it satisfies the Baum-Connes conjecture. Since $R \rtimes \Gamma$ is also torsionfree, we know by \cite[Proposition~6.3.1]{Val} that $\tau^{R \rtimes \Gamma}_*(K_0(C^*(R \rtimes \Gamma))) = \Zz$. So we obtain a decomposition $K_0(C^*(R \rtimes \Gamma)) = \Zz[1] \oplus \ker(\tau^{R \rtimes \Gamma}_*)$. By Lemma~\ref{lem-res3} in the appendix, we know that $\tau^{\mfp^{h_{\mfp}} \rtimes \Gamma}_* \circ \res_{R \rtimes \Gamma}^{\mfp^{h_{\mfp}} \rtimes \Gamma} = N(\mfp)^{h_{\mfp}} \cdot \tau^{R \rtimes \Gamma}_*$. Moreover, we certainly have $\tau^{R \rtimes \Gamma}_* \circ (\mu_{a(\mfp^{h_{\mfp}})})_* = \tau^{\mfp^{h_{\mfp}} \rtimes \Gamma}_*$. So we conclude that $\tau^{R \rtimes \Gamma}_* \circ (\mu_{a(\mfp^{h_{\mfp}})})_* \circ \res_{R \rtimes \Gamma}^{\mfp^{h_{\mfp}} \rtimes \Gamma} = N(\mfp)^{h_{\mfp}} \cdot \tau^{R \rtimes \Gamma}_*$. This shows that $(\mu_{a(\mfp^{h_{\mfp}})})_* \circ \res_{R \rtimes \Gamma}^{\mfp^{h_{\mfp}} \rtimes \Gamma} (\ker(\tau^{R \rtimes \Gamma}_*)) \subseteq \ker(\tau^{R \rtimes \Gamma}_*)$. Therefore, $(\id - (\mu_{a(\mfp^{h_{\mfp}})})_* \circ \res_{R \rtimes \Gamma}^{\mfp^{h_{\mfp}} \rtimes \Gamma})$ preserves the direct sum decomposition $K_0(C^*(R \rtimes \Gamma)) = \Zz[1] \oplus \ker(\tau^{R \rtimes \Gamma}_*)$. It follows that
\bglnoz
  && (\id - (\mu_{a(\mfp^{h_{\mfp}})})_* \circ \res_{R \rtimes \Gamma}^{\mfp^{h_{\mfp}} \rtimes \Gamma})(K_0(C^*(R \rtimes \Gamma))) \cap \Zz[1] \\
  &=& \Zz(\id - (\mu_{a(\mfp^{h_{\mfp}})})_* \circ \res_{R \rtimes \Gamma}^{\mfp^{h_{\mfp}} \rtimes \Gamma})([1]) = ((N(\mfp)^{h_{\mfp}}-1)\Zz)[1].
\eglnoz
\eproof

Let $m$ be the number of roots of unity in $K$, and let $1$ now denote the unit of $C^*(R \rtimes R^*)$.
\blemma
\label{tors-max}
\bgloz
  (\id - (\mu_{a(\mfp^{h_{\mfp}})})_* \circ \res_{R \rtimes R^*}^{\mfp^{h_{\mfp}} \rtimes R^*})(K_0(C^*(R \rtimes R^*))) \cap \Zz[1]
  \subseteq \rukl{\tfrac{N(\mfp)^{h_{\mfp}} - 1}{\gcd(m,N(\mfp)^{h_{\mfp}} - 1)} \cdot \Zz} [1].
\egloz
\elemma
\bproof
By Lemma~\ref{lem-res2} in the appendix, we know that
\bgloz
  \res_{R \rtimes R^*}^{R \rtimes \Gamma} \circ (\mu_{a(\mfp^{h_{\mfp}})})_* 
  = (\mu^{\Gamma}_{a(\mfp^{h_{\mfp}})})_* \circ \res_{\mfp^{h_{\mfp}} \rtimes R^*}^{\mfp^{h_{\mfp}} \rtimes \Gamma}.
\egloz
In addition, it is clear that
\bgloz
  \res_{\mfp^{h_{\mfp}} \rtimes R^*}^{\mfp^{h_{\mfp}} \rtimes \Gamma} \circ \res_{R \rtimes R^*}^{\mfp^{h_{\mfp}} \rtimes R^*}
  = \res_{R \rtimes R^*}^{\mfp^{h_{\mfp}} \rtimes \Gamma}
  = \res_{R \rtimes \Gamma}^{\mfp^{h_{\mfp}} \rtimes \Gamma} \circ \res_{R \rtimes R^*}^{R \rtimes \Gamma}.
\egloz
Therefore, the following diagram commutes:
\bgloz
  \xymatrix@C=1mm{
  K_0(C^*(R \rtimes R^*)) 
  \ar[d]_{\res_{R \rtimes R^*}^{R \rtimes \Gamma}} 
  \ar[rr]^{\id - (\mu_{a(\mfp^{h_{\mfp}})})_* \circ \res_{R \rtimes R^*}^{\mfp^{h_{\mfp}} \rtimes R^*}} 
  & \ \ \ \ \ \ \ \ \ \ \ \ & K_0(C^*(R \rtimes R^*)) 
  \ar[d]^{\res_{R \rtimes R^*}^{R \rtimes \Gamma}} \\
  K_0(C^*(R \rtimes \Gamma)) 
  \ar[rr]^{\id - (\mu^{\Gamma}_{a(\mfp^{h_{\mfp}})})_* \circ \res_{R \rtimes \Gamma}^{\mfp^{h_{\mfp}} \rtimes \Gamma}} 
  & \ \ \ \ \ \ \ \ \ \ \ \ & K_0(C^*(R \rtimes \Gamma))
  }
\egloz
Now take $x \in K_0(C^*(R \rtimes R^*))$ with $(\id - (\mu_{a(\mfp^{h_{\mfp}})})_* \circ \res_{R \rtimes R^*}^{\mfp^{h_{\mfp}} \rtimes R^*})(x) = z[1_{C^*(R \rtimes R^*)}] \in \Zz[1_{C^*(R \rtimes R^*)}]$. Then
\bglnoz
  m z [1_{C^*(R \rtimes \Gamma)}] &=& \res_{R \rtimes R^*}^{R \rtimes \Gamma} (z [1_{C^*(R \rtimes R^*)}]) \\
  &=& \res_{R \rtimes R^*}^{R \rtimes \Gamma} ((\id - (\mu_{a(\mfp^{h_{\mfp}})})_* \circ \res_{R \rtimes R^*}^{\mfp^{h_{\mfp}} \rtimes R^*})(x)) \\
  &=& (\id - (\mu^{\Gamma}_{a(\mfp^{h_{\mfp}})})_* \circ \res_{R \rtimes \Gamma}^{\mfp^{h_{\mfp}} \rtimes \Gamma}) (\res_{R \rtimes R^*}^{R \rtimes \Gamma}(x)) \\
  &\in& (\id - (\mu^{\Gamma}_{a(\mfp^{h_{\mfp}})})_* \circ \res_{R \rtimes \Gamma}^{\mfp^{h_{\mfp}} \rtimes \Gamma}) K_0(C^*(R \rtimes \Gamma)) 
  \cap \Zz[1_{C^*(R \rtimes \Gamma)}] \\
  &=& ((N(\mfp)^{h_{\mfp}} - 1)\Zz)[1_{C^*(R \rtimes \Gamma)}]
\eglnoz
by the previous lemma.

Hence $mz \in (N(\mfp)^{h_{\mfp}} - 1)\Zz$ which implies $z \in \frac{N(\mfp)^{h_{\mfp}} - 1}{\gcd(m,N(\mfp)^{h_{\mfp}} - 1)} \cdot \Zz$.
\eproof

Let $\zeta$ be a root of unity which generates $\mu$, i.e. $\mu = \spkl{\zeta}$. Let $\mfp$ be a prime ideal with the property that for every $1 \leq i \leq m-1$, we have $1 - \zeta^i \notin \mfp$. This implies that the order of $\zeta$ in $R / \mfp$ is $m$. Thus $(R / \mfp)^*$ contains a cyclic subgroup of order $m$. Therefore $m$ divides $N(\mfp) - 1$, hence also $N(\mfp)^{h_{\mfp}} - 1$. This means that $\gcd(m,N(\mfp)^{h_{\mfp}} - 1) = m$.
\blemma
For a prime ideal $\mfp$ with $1 - \zeta^i \notin \mfp$ for every $1 \leq i \leq m-1$, we have
\bgloz
  \tfrac{N(\mfp)^{h_{\mfp}} - 1}{m} [1] \in (\id - (\mu_{a(\mfp^{h_{\mfp}})})_* \circ \res_{R \rtimes R^*}^{\mfp^{h_{\mfp}} \rtimes R^*})(K_0(C^*(R \rtimes R^*))).
\egloz
Here $1$ is the unit in $C^*(R \rtimes R^*)$.
\elemma
\bproof
Let $e$ be the projection $\tfrac{1}{m} \sum_{i=0}^{m-1} (u_{(0,\zeta)})^i$ in $C^*(R \rtimes R^*)$. Consider the homomorphism $\varphi: C^*(R \rtimes R^*) \to M_{N(\mfp)}(C^*(\mfp^{h_{\mfp}} \rtimes R^*))$ for $G = R \rtimes R^*$, $H = \mfp^{h_{\mfp}} \rtimes R^*$ from \eqref{phi-res} in the appendix. Now $\varphi(u_{(0,\zeta)})$ is a matrix which can be decomposed into irreducible ones. In this decomposition, we only obtain one $1$-dimensional (i.e. $1 \times 1$) irreducible matrix, namely corresponding to the trivial coset $\mfp^{h_{\mfp}} \rtimes R^*$ of $R \rtimes R^* / \mfp^{h_{\mfp}} \rtimes R^*$, and $\tfrac{N(\mfp)^{h_{\mfp}} - 1}{m}$ $m$-dimensional (i.e. $m \times m$) irreducible matrices. This follows from our assumption that $1 - \zeta^i \notin \mfp$ for all $1 \leq i \leq m-1$. The analogues of Lemma~5.6 and Lemma~5.7 in \cite{Li-Lu} imply that
\bgloz
  ((\mu_{a(\mfp^{h_{\mfp}})})_* \circ \res_{R \rtimes R^*}^{\mfp^{h_{\mfp}} \rtimes R^*}) [e] = [e] + \tfrac{N(\mfp)^{h_{\mfp}} - 1}{m} [1].
\egloz
\eproof
As an immediate consequence of the two previous lemmas, we obtain
\bcor
\label{i-mrcap1}
For every prime ideal $\mfp$ with $1 - \zeta^i \notin \mfp$ for all $1 \leq i \leq m-1$,
\bgloz
  (\id - (\mu_{a(\mfp^{h_{\mfp}})})_* \circ \res_{R \rtimes R^*}^{\mfp^{h_{\mfp}} \rtimes R^*})(K_0(C^*(R \rtimes R^*))) \cap \Zz[1]
  = (\tfrac{N(\mfp)^{h_{\mfp}} - 1}{m} \cdot \Zz) [1].
\egloz
\ecor
\bprop
\label{tors}
For every prime ideal $\mfp$ with $1 - \zeta^i \notin \mfp$ for all $1 \leq i \leq m-1$, the $K_0$-class $[1]$ of the unit of $C^*_r(R \rtimes R\reg) / I_{\gekl{\mfp}}$ is a torsion element of order $\tfrac{N(\mfp)^{h_{\mfp}} - 1}{m}$.
\eprop
\bproof
Using the six term exact sequence in K-theory for the short exact sequence $0 \to I_{\gekl{\mfp}} \overset{i_{\gekl{\mfp}}}{\lori} C^*_r(R \rtimes R\reg) \to C^*_r(R \rtimes R\reg) / I_{\gekl{\mfp}} \to 0$, our claim follows once we know that
\bgl
\label{imcap1}
  \img(i_{\gekl{\mfp}}) \cap \Zz[1_{C^*_r(R \rtimes R\reg)}] = (\tfrac{N(\mfp)^{h_{\mfp}} - 1}{m} \cdot \Zz)[1_{C^*_r(R \rtimes R\reg)}]
\egl
in $K_0(C^*_r(R \rtimes R\reg))$.

As before, let $[1]_{[R]}$ be the element of $\bigoplus_{\mfk \in Cl_K} K_0(C^*(\mfa_{\mfk} \rtimes R^*))$ defined by
\bgloz
  ([1]_{[R]})_{\mfk} = 
  \bfa
    0 \falls \mfk \neq [R], \\
    [1_{C^*(R \rtimes R^*)}] \falls \mfk = [R].
  \efa
\egloz
Using the identifications
\bgloz
  \sum_{\mfk} (\kappa_{\mfk})_*: \bigoplus_{\mfk \in Cl_K} K_*(C^*(\mfa_{\mfk} \rtimes R^*)) \cong K_*(I_{\gekl{\mfp}}),
\egloz
\bgloz
  \sum_{\mfk} (\iota_{\mfk})_*: \bigoplus_{\mfk \in Cl_K} K_*(C^*(\mfa_{\mfk} \rtimes R^*)) \cong K_*(C^*_r(R \rtimes R\reg))
\egloz
and the fact $(\sum_{\mfk} (\iota_{\mfk})_*)[1]_{[R]} = [1_{C^*_r(R \rtimes R\reg)}]$, \eqref{imcap1} is obviously equivalent to
\bgloz
  \img \rukl{(\sum_{\ti{\mfk}} (\iota_{\ti{\mfk}})_*)^{-1} \circ i_{\gekl{\mfp}} \circ (\sum_{\mfk} (\kappa_{\mfk})_*)} \cap \Zz [1]_{[R]} 
  = (\tfrac{N(\mfp)^{h_{\mfp}} - 1}{m} \cdot \Zz) [1]_{[R]}.
\egloz
But this is precisely what we obtain by combining Lemma~\ref{imim} with Corollary~\ref{i-mrcap1}.
\eproof

Set $p_{\max} \defeq \max \menge{p \in \Nz \text{ prime}}{p \vert N(1 - \zeta^i) \text{ for some } 1 \leq i \leq m-1}$. Also let $n \defeq [K:\Qz]$ and let $h$ be the class number of $K$. Moreover, let $\Prim_{\min}$ be the set of minimal non-zero primitive ideals of $C^*_r(R \rtimes R\reg)$. The following result tells us how to read off splitting numbers from K-theoretic invariants, at least for all but finitely many prime numbers.
\blemma
For every prime $p \in \Nz$ with $\tfrac{p-1}{m} > p_{\max}^{nh} - 1$, we have the following formula for the splitting number $g_K(p)$ of $p$ in $K$:
\bgloz
  g_K(p) = \# \menge{I \in \Prim_{\min}}{p \vert (m \cdot \ord(1_{C^*_r(R \rtimes R\reg) / I}) + 1)}.
\egloz
\elemma
\bproof
To prove \an{$\leq$}, observe that $\tfrac{p-1}{m} > p_{\max}^{nh} - 1$ implies that $p-1 > p_{\max} - 1$, hence $p > p_{\max}$. Now assume that a prime ideal $\mfp$ lies above $p$, i.e. $\mfp \cap \Zz = p \Zz$. Then we must have $1 - \zeta^i \notin \mfp$ for all $1 \leq i \leq m-1$ since otherwise, $1 - \zeta^i \in \mfp$ would imply $(1 - \zeta^i)R \subseteq \mfp$, hence $N(\mfp) \vert N(1 - \zeta^i)$, thus $p \vert N(1 - \zeta^i)$ for some $1 \leq i \leq m-1$. But this contradicts $p > p_{\max}$. Thus every $\mfp \in \cP$ with $\mfp \cap \Zz = p \Zz$ satisfies $\ord(1_{C^*_r(R \rtimes R\reg) / I_{\gekl{\mfp}}}) = \tfrac{N(\mfp)^{h_{\mfp}} - 1}{m}$. Hence $p$ divides $m \cdot \ord(1_{C^*_r(R \rtimes R\reg) / I_{\gekl{\mfp}}}) + 1$. This shows \an{$\leq$}.

To see \an{$\geq$}, let $\mfp \in \cP$ be a prime ideal such that $p$ divides $m \cdot \ord(1_{C^*_r(R \rtimes R\reg) / I_{\gekl{\mfp}}}) + 1$. This implies $m \cdot \ord(1_{C^*_r(R \rtimes R\reg) / I_{\gekl{\mfp}}}) + 1 \geq p$, hence $\ord(1_{C^*_r(R \rtimes R\reg) / I_{\gekl{\mfp}}}) \geq \tfrac{p-1}{m} > p_{\max}^{nh} - 1$. Let $q \in \Nz$ be the prime number determined by $\mfp \cap \Zz = q \Zz$. Then we know by Lemma~\ref{tors-max} that $q^{nh} - 1 \geq \ord(1_{C^*_r(R \rtimes R\reg) / I_{\gekl{\mfp}}}) \geq \tfrac{p-1}{m} > p_{\max}^{nh} - 1$. This implies $q > p_{\max}$. Hence we deduce that $1 - \zeta^i \notin \mfp$ for all $1 \leq i \leq m-1$, because otherwise, we would again obtain $N(\mfp) \vert N(1 - \zeta^i)$ and hence $q \vert N(1 - \zeta^i)$ for some $1 \leq i \leq m-1$. But this cannot happen since $q > p_{\max}$. Hence $\ord(1_{C^*_r(R \rtimes R\reg) / I_{\gekl{\mfp}}}) = \tfrac{N(\mfp)^{h_{\mfp}} - 1}{m}$. Therefore, our assumption that $p$ divides $m \cdot \ord(1_{C^*_r(R \rtimes R\reg) / I_{\gekl{\mfp}}}) + 1$ implies that $p$ divides $N(\mfp)^{h_{\mfp}}$, hence $\mfp \cap \Zz = p \Zz$ (and thus $p=q$). This shows \an{$\geq$}.
\eproof

As a consequence, we obtain our main result Theorem~\ref{main}.
\btheooz[Theorem~\ref{main}]
Assume that $K$ and $L$ are two number fields which have the same number of roots of unity. Let $R$ and $S$ be their rings of integers. If $C^*_r(R \rtimes R\reg) \cong C^*_r(S \rtimes S\reg)$, then $K$ and $L$ are arithmetically equivalent.
\etheooz
\bproof
It follows from our assumptions and the previous lemma that for all but finitely many prime numbers $p \in \Nz$, we have $g_K(p) = g_L(p)$. Hence $K$ and $L$ are arithmetically equivalent by \cite[Main Theorem]{Per-Stu}.
\eproof

And we also obtain Theorem~\ref{mainGal}.
\btheooz[Theorem~\ref{mainGal}]
Assume that $K$ and $L$ are finite Galois extensions of $\Qz$ which have the same number of roots of unity. Then we have $C^*_r(R \rtimes R\reg) \cong C^*_r(S \rtimes S\reg)$ if and only if $K \cong L$.
\etheooz
\bproof
The direction \an{$\Larr$} is clear: If $K \cong L$, then $R \cong S$ and thus $R \rtimes R\reg \cong S \rtimes S\reg$. This implies $C^*_r(R \rtimes R\reg) \cong C^*_r(S \rtimes S\reg)$. The other direction \an{$\Rarr$} follows from Theorem~\ref{main} and the fact that two finite Galois extensions of $\Qz$ which are arithmetically equivalent must already be isomorphic (see \cite{Per} or \cite[Chapter~VII, Corollary~(13.10)]{Neu}).
\eproof

\section{The number of roots of unity}
\label{number-rou}

An obvious question is whether in the same setting as in Theorem~\ref{main} and Theorem~\ref{mainGal}, we can deduce from $C^*_r(R \rtimes R\reg) \cong C^*_r(S \rtimes S\reg)$ that $K$ and $L$ have the same number of roots of unity. If this was the case, then we could drop this extra assumption in our main result. Unfortunately, we cannot answer this question in general. However, in the following, we provide a positive answer in a special case, namely for purely imaginary number fields.

We start with the following easy observation, which is an immediate consequence of our results from the previous section.

\blemma
Let $K$ be a number field with ring of integers $R$. Then we have the following formula for the degree of $K$ over $\Qz$:
\bgloz
  [K:\Qz] = \limsup_{T \to \infty} \# \menge{I \in \Prim_{\min}}{\ord(1_{C^*_r(R \rtimes R\reg) / I}) = T}.
\egloz
\elemma
\bproof
By \cite[Chapter~I, \S~8]{Neu}, we know that there are infinitely many prime numbers which totally split in $K$. This observation, combined with Proposition~\ref{tors}, gives the desired formula for the degree.
\eproof

Let $K$ be a number field with ring of integers $R$. In the previous section, we only studied the minimal non-zero primitive ideals of $C^*_r(R \rtimes R\reg)$. Now we turn to the maximal primitive ideal. As observed in \cite[Remark~3.9]{Ech-La}, the corresponding quotient is canonically isomorphic to the ring C*-algebra $\fA[R]$ from \cite{Cu-Li1}. Let $\pi: C^*_r(R \rtimes R\reg) \onto \fA[R]$ be the canonical projection, and let $\pi_*$ be the induced homomorphism on $K_0$.
\bprop
Let $K$ be a purely imaginary number field, i.e., $K$ has no real embeddings. Let $m$ be the number of roots of unity in $K$. Then we can compute the rank of $\pi_*$ as follows:
\bgloz
  \rk(\pi_*) = m 2^{\halb [K:\Qz] - 2}.
\egloz
\eprop
\bproof
The proof consists of two steps. First, consider the identification $\fA[R] \cong C(\Rbar) \rtimes^e (R \rtimes R\reg)$ from \cite{Cu-Li1}. Here $\Rbar$ is the profinite completion of $R$, $\rtimes^e$ stands for \an{semigroup crossed product by endomorphisms}, and $R \rtimes R\reg$ acts on $C(\Rbar)$ by affine transformations. Using this identification, we obtain a canonical homomorphism $C(\Rbar) \rtimes R \rtimes R^* \to \fA[R]$, hence a homomorphism $K_0(C(\Rbar) \rtimes R \rtimes R^*) \to K_0(\fA[R])$. Let $r$ be the rank of the image of this homomorphism. The first step of our argument is to show $\rk(\pi_*) = r$.

Here is why we have $\rk(\pi_*) = r$: As before, given $\mfk \in Cl_K$, let $\mfa_{\mfk}$ be an ideal of $R$ which represents $\mfk$. Then we have for every $\mfk \in Cl_K$ a commutative diagram
\bgl
\label{CD'}
  \xymatrix@C=1mm{
  K_0(C^*_r(\mfa_{\mfk} \rtimes R^*)) \ar[rr] \ar[d] & \ \ \ \ \ \ & K_0(C(\Rbar) \rtimes R \rtimes R^*) \ar[d] \\
  K_0(C^*_r(R \rtimes R\reg)) \ar[rr] & \ \ \ \ \ \ & K_0(\fA[R])
  }
\egl
All the arrows are induced by canonical homomorphisms on the level of C*-algebras. Now Theorem~8.2.1 in \cite{C-E-L1} tells us that the images of the left vertical arrows, for all $\mfk \in Cl_K$, sum up to $K_0(C^*_r(R \rtimes R\reg))$. This fact, together with commutativity of the diagram, immediately yields the inequality $\rk(\pi_*) \leq r$. For the reverse inequality, note that the analysis in \cite[\S~6]{Cu-Li2} and \cite[\S~4]{Li-Lu} show that $K_0(C(\Rbar) \rtimes R \rtimes R^*)$ can be written as a stationary inductive limit with groups given by $K_0(C^*(R \rtimes R^*))$ and structure maps of a particular form. It follows from this description that the canonical homomorphism $C^*(R \rtimes R^*) \to C(\Rbar) \rtimes R \rtimes R^*$ induces a rationally surjective homomorphism $K_0(C^*(R \rtimes R^*)) \to K_0(C(\Rbar) \rtimes R \rtimes R^*)$ in K-theory. All we have to do now is to consider the commutative diagram \eqref{CD'} for $\mfk = [R]$, i.e., $\mfa_{\mfk} = R$. Since the upper horizontal arrow is rationally surjective, we immediately obtain $\rk(\pi_*) \geq r$, as desired.

Secondly, the duality theorem from \cite[\S~4]{Cu-Li1} allows us to identify $K_0(C(\Rbar) \rtimes R \rtimes R^*)$ with $K_0(C_0(\Az_{\infty}) \rtimes K \rtimes R^*)$ and $K_0(\fA[R])$ with $K_0(C_0(\Az_{\infty}) \rtimes K \rtimes K\reg)$ such that the diagram
\bgloz
  \xymatrix@C=1mm{
  K_0(C(\Rbar) \rtimes R \rtimes R^*) \ar[rr] \ar[d]^{\cong} & \ \ \ \ \ \ & K_0(\fA[R]) \ar[d]^{\cong} \\
  K_0(C_0(\Az_{\infty}) \rtimes K \rtimes R^*) \ar[rr] & \ \ \ \ \ \ & K_0(C_0(\Az_{\infty}) \rtimes K \rtimes K\reg)
  }
\egloz
commutes. Here $\Az_{\infty}$ is the infinite adele space over $K$. This commutative diagram shows that $r$ coincides with the rank of the image of $K_0(C_0(\Az_{\infty}) \rtimes K \rtimes R^*) \to K_0(C_0(\Az_{\infty}) \rtimes K \rtimes K\reg)$. By \cite[Corollary~4.15]{Li-Lu}, we know that for every $c \in \Zz\reg \subseteq R\reg$, $c>1$, the canonical homomorphism $K_0(C_0(\Az_{\infty}) \rtimes K \rtimes (R^* \times \spkl{c})) \to K_0(C_0(\Az_{\infty}) \rtimes K \rtimes K\reg)$ is rationally injective. Hence $r$ also concides with the rank of the image of $K_0(C_0(\Az_{\infty}) \rtimes K \rtimes R^*) \to K_0(C_0(\Az_{\infty}) \rtimes K \rtimes (R^* \times \spkl{c}))$. Using this, we now show that $r = m 2^{\halb [K:\Qz] - 2}$.

Consider the commutative diagram
\bgloz
  \xymatrix@C=1mm{
  K_0(C_0(\Az_{\infty}) \rtimes K \rtimes R^*) \ar[rr] & \ \ \ \ \ \ & K_0(C_0(\Az_{\infty}) \rtimes K \rtimes (R^* \times \spkl{c})) \\
  K_0(C_0(\Az_{\infty}) \rtimes R^*) \ar[rr] \ar[u] & \ \ \ \ \ \ & K_0(C_0(\Az_{\infty}) \rtimes (R^* \times \spkl{c})) \ar[u]
  }
\egloz
Again, all the arrows are induced by canonical homomorphisms on the level of C*-algebras. Using the Pimsner-Voiculescu exact sequence, it is easy to see that the horizontal arrow at the bottom of the diagram is rationally injective. By the injectivity results in \cite[\S~5]{Li-Lu}, we also know that the right vertical arrow in the diagram is rationally injective. This shows that $r \geq \rk(K_0(C_0(\Az_{\infty}) \rtimes R^*)) = m 2^{\halb [K:\Qz] - 2}$. For the last equality, we have used the assumption that $K$ is purely imaginary, so that $\rk(R^*) = \halb [K:\Qz] - 1$. At the same time, we know that $\rk(K_0(C_0(\Az_{\infty}) \rtimes K \rtimes (R^* \times \spkl{c}))) = m 2^{\halb [K:\Qz] - 1}$ by \cite[Corollary~4.15]{Li-Lu}. And finally, examining the Pimsner-Voiculescu exact sequence for $C_0(\Az_{\infty}) \rtimes K \rtimes (R^* \times \spkl{c}) \cong (C_0(\Az_{\infty}) \rtimes K \rtimes R^*) \rtimes \Zz$, we deduce that $r \leq m 2^{\halb [K:\Qz] - 2}$.
\eproof

\bcor
Let $K$ and $L$ be two number fields with rings of integers $R$ and $S$, respectively. If $C^*_r(R \rtimes R\reg) \cong C^*_r(S \rtimes S\reg)$, then $K$ and $L$ must have the same degree over $\Qz$. Moreover, if $K$ and $L$ are both purely imaginary, then $K$ and $L$ have the same number of roots of unity.
\ecor
\bproof
The first part of our assertion is an immediate consequence of the previous lemma. For the second claim, note that an isomorphism $C^*_r(R \rtimes R\reg) \cong C^*_r(S \rtimes S\reg)$ must preserve the maximal primitive ideals. Therefore, the previous proposition implies that $K$ and $L$ must have the same number of roots of unity.
\eproof

With the help of this corollary, we obtain the following stronger version of our main result:
\btheooz[Theorem~\ref{stronger1}]
Let $K$ and $L$ be finite Galois extensions of $\Qz$ with rings of integers $R$ and $S$, respectively. Assume that either both $K$ and $L$ have at least one real embedding, or that both $K$ and $L$ are purely imaginary. Then $C^*_r(R \rtimes R\reg) \cong C^*_r(S \rtimes S\reg)$ if and only if $K \cong L$.
\etheooz

\bremark
In the last theorem, we can replace the assumption that either both $K$ and $L$ have at least one real embedding or that both $K$ and $L$ are purely imaginary by the assumption that either both $K$ and $L$ have only two roots of unity $+1$ and $-1$, or that $K$ and $L$ both have more than two roots of unity. The reason is that a number field with more than two roots of unity can never have a real embedding. In this way, we obtain Theorem~\ref{stronger2}.

However, the question remains whether $C^*_r(R \rtimes R\reg) \cong C^*_r(S \rtimes S\reg)$ excludes the possibility that $K$ has roots of unity $+1$ and $-1$ whereas $L$ has more than two roots of unity.
\eremark

\begin{appendix}

\section{Restriction homomorphisms}

In the following, we collect a few certainly well-known facts about restriction homomorphisms in $K_0$. Of course, we do not claim any originality here.

Assume that $G$ is a group and that $H$ is a subgroup of $G$ of finite index. Let $N \defeq (G:H)$. The restriction $\res_G^H: K_0(C^*_r(G)) \to K_0(C^*_r(H))$ is defined by sending a finitely generated projective right $C^*_r(G)$-module $M$ to $M \vert_{C^*_r(H)}$, i.e. we just restrict the $C^*_r(G)$-action on $M$ to a $C^*_r(H)$-action. Let us now describe $\res_G^H$ in terms of homomorphisms of C*-algebras. Let $\cR \subseteq G$ be a full set of representatives of $G / H$, i.e $\cR \ni \gamma \ma \gamma H \in G/H$ is a bijection. From the set-theoretic bijection $G = \dotcup_{\gamma \in \cR} \gamma H \cong \dotcup_{\gamma \in \cR} H$, we obtain a unitary
\bgl
\label{unitary}
  \ell^2(G) \cong \bigoplus_{\gamma \in \cR} \ell^2(H), \ \sum_{\gamma \in \cR, h \in H} \lambda_{\gamma h} \ve_{\gamma h} \ma (\sum_h \lambda_{\gamma h} \ve_h)_{\gamma}.
\egl
Let $E_H$ be the orthogonal projection in $\cL(\ell^2(G))$ onto $\ell^2(H) \subseteq \ell^2(G)$. Conjugation with the unitary from \eqref{unitary} yields an isomorphism
\bgl
\label{identification}
  \cL(\ell^2(g)) \cong M_N(\cL(\ell^2(H))), \ T \ma (E_H \lambda_{\gamma}^{-1} T \lambda_{\gamma'} E_H)_{\gamma, \gamma'}.
\egl
Let us now represent $C(G/H) \rtimes_r G$ in a canonical and faithful way on $\ell^2(G)$ and identify $C(G/H) \rtimes_r G$ with its image in $\cL(\ell^2(G))$. Then the isomorphism from \eqref{identification} identifies $C(G/H) \rtimes_r G$ with $M_N(C^*_r(H)) \subseteq M_N(\cL(\ell^2(H)))$. Let us call this identification $\phi$, and let $\varphi$ be the composition
\bgl
\label{phi-res}
  \varphi: C^*_r(G) \to C(G/H) \rtimes_r G \underset{\cong}{\overset{\phi}{\lori}} M_N(C^*_r(H)).
\egl
Moreover, let $e$ be the minimal projection in $M_N(C^*_r(H))$ corresponding to the entry in the upper left corner, so that $e \otimes \id: C^*_r(H) \to M_N(C^*_r(H))$ is the canonical embedding into the upper left corner. We then have the following
\blemma
\label{lem-res1}
$\res_G^H = (e \otimes \id)_*^{-1} \circ \varphi_*$.
\elemma
\bproof
We have to show that for every finitely generated projective right $C^*_r(G)$-module $M$, we can identify
\bgloz
  M \vert_{C^*_r(H)} \otimes_{e \otimes \id} M_N(C^*_r(H)) \text{ and } M \otimes_{\varphi} M_N(C^*_r(H))
\egloz
as right $M_N(C^*_r(H))$-modules.

But it is easy to see that the homomorphisms
\bgloz
  M \vert_{C^*_r(H)} \otimes_{e \otimes \id} M_N(C^*_r(H)) \to M \otimes_{\varphi} M_N(C^*_r(H)), m \otimes x \ma m \otimes ex
\egloz
and
\bgloz
  M \otimes_{\varphi} M_N(C^*_r(H)) \to M \vert_{C^*_r(H)} \otimes_{e \otimes \id} M_N(C^*_r(H)), 
  m \otimes x \ma \sum_{\gamma \in \cR} m \lambda_{\gamma} \otimes \varphi(\lambda_{\gamma^{-1}}) x
\egloz
are mutually inverse homomorphisms of right $M_N(C^*_r(H))$-modules.
\eproof

Now let $H \subseteq G$ and $H' \subseteq G'$ be two subgroups of groups such that $(G:H), (G':H') < \infty$. Assume that $\alpha: G \to G'$ is a group isomorphism with $\alpha(H) = H'$. Such an isomorphism induces homomorphisms $C^*_r(\alpha): C^*_r(G) \to C^*_r(G')$ and $C^*_r(\alpha \vert_H): C^*_r(H) \to C^*_r(H')$ as well as homomorphisms $\alpha_* = (C^*_r(\alpha))_*: K_0(C^*_r(G)) \to K_0(C^*_r(G'))$, $(\alpha \vert_H)_* = (C^*_r(\alpha \vert H))_*: K_0(C^*_r(H)) \to K_0(C^*_r(H'))$.
\blemma
\label{lem-res2}
$(\alpha \vert_H)_* \circ \res_G^H = \res_{G'}^{H'} \circ \alpha_*$.
\elemma
\bproof
Given a right $C^*_r(G)$-module $M$, we can obviously identify $M \vert_{C^*_r(H)} \otimes_{C^*_r(\alpha \vert_H)} C^*_r(H')$ with $(M \otimes_{C^*_r(\alpha)} C^*_r(G')) \vert_{C^*_r(H')}$ by the mutually inverse homomorphisms $m \otimes x \ma m \otimes x$ and $m C^*_r(\alpha)^{-1}(y) \otimes 1 \mafr m \otimes y$.
\eproof

Finally, we need the following observation about traces: Let $\tau^G$ and $\tau^H$ denote the canonical tracial states on $C^*_r(G)$ and $C^*_r(H)$. We denote the induced homomorphisms on $K_0$ by $\tau^G_*$ and $\tau^H_*$.
\blemma
\label{lem-res3}
$\tau^H_* \circ \res_G^H = N \cdot \tau^G_*$.
\elemma
\bproof
Let $M_N(\tau^H)$ be the trace on $M_N(C^*_r(H))$ such that $M_N(\tau^H) \circ (e \otimes \id) = \tau^H$. An easy computation shows that $M_N(\tau^H) \circ \varphi = N \cdot \tau^G$, where $\varphi$ is defined in \eqref{phi-res}. Our claim follows.
\eproof

\end{appendix}

\end{document}